\documentclass[11pt,twoside,reqno]{amsart}
\numberwithin{equation}{section}
\usepackage[centertags]{amsmath}
\usepackage{amssymb,amsfonts,amsthm}
\IfFileExists{mathrsfs.sty}{
\usepackage{mathrsfs}}         
   {
     \IfFileExists{eucal.sty}{
        \usepackage[mathscr]{eucal}} 
   {
   }
}
 
\newcommand{\cA}{{\mathcal A}}
\newcommand{\cB}{{\mathcal B}}
\newcommand{\cJ}{{\mathcal J}}
\newcommand{\cE}{{\mathscr E}}
\newcommand{\cC}{{\mathcal C}}
\newcommand{\cD}{{\mathcal D}}
\newcommand{\cQ}{{\mathcal Q}}
\newcommand{\cP}{{\mathcal P}}
\newcommand{\cS}{{\mathscr S}}

\newcommand{\N}{\mathbb{N}}
\newcommand{\Z}{\mathbb{Z}}
\newcommand{\R}{\mathbb{R}}
\newcommand{\C}{\mathbb{C}}
\newcommand{\GL}{\operatorname{GL}}

\newcommand{\End}{\operatorname{End}}

\newcommand{\character}{\operatorname{char}} 

\newcommand{\id}{\operatorname{Id}}
\newcommand{\im}{\operatorname{Im}}
\newcommand{\Ker}{\operatorname{Ker}}

\newcommand{\tot}{\operatorname{Tot}}
\newcommand{\tr}{\operatorname{tr}}

\newcommand{\TR}{\operatorname{TR}}
\newcommand{\fTR}{\overline{\operatorname{TR}}}
\newcommand{\lTR}{\widetilde{\operatorname{TR}}}
\newcommand{\ch}{\operatorname{ch}}
\newcommand{\slch}{\,/\!\!\!\!\operatorname{ch}}
\newcommand{\sllch}{\,/\!\!/\!\!\!\!\!\operatorname{ch}}
\newcommand{\dbar}{{d \! \bar{ } \,}}
  
\newcommand{\DF}{\operatorname{DF}}
\newcommand{\SF}{\operatorname{SF}}
\newcommand{\Ell}{\operatorname{Ell}}
\newcommand{\Proj}{\operatorname{P}}
\newcommand{\AP}{\operatorname{AP}} 
\newcommand{\approj}{(\AP_\infty(\cA),\Proj_\infty(\cA))}

\newcommand{\CS}{\operatorname{CS}}

\newcommand{\sym}{\operatorname{S}}
\newcommand{\PS}{\operatorname{PS}}
\newcommand{\polyn}{\C[\mu_1,\ldots,\mu_p]}
\newcommand{\csym}{\textup{C}\mathcal S}
\newcommand{\symb}{\operatorname{Symb}}
\newcommand{\CL}{\textup{CL}} 
\newcommand{\pdo}{\textup{L}} 
\newcommand{\Op}{\operatorname{Op}}
\newcommand{\regint}{-\hspace*{-1em}\int}

\newcommand{\LIM}{\operatorname*{LIM}}
\newcommand{\romlabel}{\renewcommand{\labelenumi}{\textup{(\roman{enumi})}}}
\newcommand{\plref}[1]{\textup{\ref{#1}}}

\theoremstyle{plain}
        \newtheorem{theorem}{Theorem}[section]
        \newtheorem{lemma}[theorem]{Lemma}
        \newtheorem{proposition}[theorem]{Proposition}

\theoremstyle{definition}
        \newtheorem{definition}[theorem]{Definition}
        \newtheorem{remark}[theorem]{Remark}
        \newtheorem{example}[theorem]{Example}

\begin{document}
\title{Relative pairing in cyclic cohomology and divisor flows}
\author{Matthias Lesch}
\thanks{The first named author was partially supported 
by \uppercase{S}onderforschungs\-be\-reich/\uppercase{T}ransregio 12 
``\uppercase{S}ymmetries and \uppercase{U}niversality in
\uppercase{M}esoscopic \uppercase{S}ystems" 
(\uppercase{B}ochum--\uppercase{D}uisburg/\uppercase{E}ssen--\uppercase{K}\"oln--\uppercase{W}arszawa)}

\address{Mathematisches Institut,
Universit\"at Bonn,
Beringstr. 1,
53115 Bonn,
Germany}

\email{lesch@math.uni-bonn.de}
\urladdr{http://www.math.uni-bonn.de/people/lesch}
\author{Henri Moscovici}
\address{Department of Mathematics,
The Ohio State University,
Columbus, OH 43210,
USA}
\email{henri@math.ohio-state.edu}

\thanks{The work of the second named author was partially
 supported by the US National Science Foundation
    award no. DMS-0245481}
\author{Markus J. Pflaum}
\address{Fachbereich Mathematik, 
Goethe-Universit\"at, 60054 Frankfurt/Main, 
Germany} 
\email{pflaum@math.uni-frankfurt.de}
\urladdr{http://www.math.uni-frankfurt.de/$\sim$pflaum}

\thanks{The third named author was partially supported by the DFG}

\subjclass[2000]{Primary 46L80; Secondary 58J30, 58J40}

\begin{abstract} 
We construct invariants of relative $K$-theory classes
of multiparameter dependent pseudodifferential operators,
which recover and generalize Melrose's divisor flow and its 
higher odd-dimensional versions
of Lesch and Pflaum. These higher divisor flows
are obtained by means of pairing the relative $K$-theory modulo 
the symbols with the cyclic cohomological characters of
relative cycles constructed out of the regularized operator trace together 
with its symbolic boundary. Besides
giving a clear and conceptual explanation to all the essential features of 
the divisor flows, this construction allows to uncover the previously
unknown even-dimensional counterparts. Furthermore, it confers to 
the totality of these invariants a purely topological interpretation, 
that of implementing  the classical Bott periodicity isomorphisms  
in a manner compatible with the suspension isomorphisms in both $K$-theory and 
in cyclic cohomology. 
We also give a precise formulation, in terms of a natural Clifford algebraic suspension,
for the relationship between the higher divisor flows and the spectral flow.

\end{abstract}

\maketitle

\section*{Introduction}  
Cyclic cohomology of associative algebras, viewed as a noncommutative
analogue of de Rham cohomology, together with its 
pairing with $K$-theory, was shown by \textsc{Connes} \cite{Con:NDG} to provide a natural
extension of the Chern-Weil construction of
characteristic classes to the general 
framework of noncommutative geometry. In this capacity, cyclic cohomology 
has been successfully exploited to
produce invariants for $K$-theory classes in a variety of interesting situations
(see \textsc{Connes} \cite{Con:NG} for an impressive array of such
applications, that include the proof of the Novikov conjecture in the case of
Gromov's word-hyperbolic groups, cf. \textsc{Connes-Moscovici} ~\cite{ConMos:NCHG}). 

In this paper we present an application of this method to the
construction of invariants of $K$-theory classes in the relative setting, which
takes advantage of the excision property in both topological 
$K$-theory and (periodic) cyclic cohomology 
(cf.~\textsc{Wodzicki} \cite{Wod:ECHRAKT},
\textsc{Cuntz--Quillen} \cite{CunQui:EBPCC}). Namely, we construct
invariants of relative $K$-theory classes
of multiparameter dependent pseudodifferential operators,
which recover and generalize the divisor flow for suspended pseudodifferential 
operators introduced by 
\textsc{Melrose} \cite{Mel:EIF}, as well as its multiparametric versions 
defined by \textsc{Lesch--Pflaum} \cite{LesPfl:TAP}.  These 
invariants  are obtained 
by pairing relative cyclic classes determined by the 
regularization \`a la Melrose of the operator trace and
its symbolic boundary  with
the relative $K$-theory of the algebras of parametric 
pseudodifferential operators and of their symbols. By the very nature of the construction, 
the `higher divisor flows' so obtained are homotopy invariant, 
additive and assume integral values.  
Besides providing a conceptual explanation for all
their essential features, this interpretation  leads naturally to the uncovering
of the formerly `missing' 
even dimensional higher eta invariants and of their
associated higher divisor flows.

To outline the origins and
content of this article in more precise terms, a modicum of notation
will be necessary. Let $M$ be a smooth compact Riemannian manifold without boundary, and
let $E$ be a hermitian vector bundle over $M$. We denote by $\CL^m(M,E)$
the classical ($1$-step polyhomogeneous) pseudodifferential operators of 
order $m$ acting between the sections of $E$. It is well-known that the 
operator trace, which is defined  on operators of order $m <  -\dim M$, 
cannot be extended to a trace on the whole algebra 
$\CL^\infty(M,E)=\bigcup\limits_{m\in\R} \CL^m(M,E)$. In fact, 
for $M^n$ connected and $n > 1$, up to a scalar multiple there is only one
tracial functional on $\CL^\infty (M,E)$, and that functional
vanishes on pseudodifferential operators of order $m <  -\dim M$,
cf. \textsc{Wodzicki} \cite{Wod:LIOSA}.
This picture changes however if one passes to `pseudodifferential suspensions'
of the algebra $\CL^\infty(M,E)$. As shown by \textsc{R.B. Melrose} \cite{Mel:EIF},
for a `natural' pseudodifferential suspension $\Psi_{\rm sus}^\infty(M,E)$ of 
$\CL^\infty(M,E)$, the operator trace on operators 
of order  $m <  -\dim M -1$
can be extended by a canonical regularization procedure to a trace on the full algebra. 
Melrose has used this regularized trace to `lift' the spectrally defined $\eta$-invariant of
\textsc{Atiyah-Patodi-Singer} \cite{APS:SARGI} to an \textit{eta} homomorphism 
$\eta: K^{\rm alg}_1 (\Psi_{\rm sus}^\infty(M,E))\rightarrow \mathbb{C}$, where
$K^{\rm alg}_1$ stands for the algebraic $K_1$-theory group. Furthermore,
by means of the variation of his generalized $\eta$-invariant, he defined
the \textit{divisor flow} between two invertibles of the algebra 
$\Psi_{\rm sus}^\infty(M,E)$ that are in the same component of the set of 
elliptic elements, and showed that it enjoys properties analogous to the 
spectral flow for self-adjoint elliptic operators.

Working with a slightly modified notion of pseudodifferential suspension, 
and for an arbitrary dimension $p\in \N$ of the parameter space, 
\textsc{Lesch} and \textsc{Pflaum} \cite{LesPfl:TAP} generalized Melrose's
trace regularization to the $ p$-fold  suspended pseudodifferential algebra
$\CL^\infty (M,E; \R^p)$ of classical parameter dependent pseudodifferential
operators. They also generalized Melrose's $\eta$-invariant to odd parametric 
dimensions, defining for $p = 2k+1$ the \textit{higher $\eta$-invariant} 
$\eta_{2k+1} (A)$ of an invertible $A \in \CL^\infty (M,E; \R^p)$. 
The appellative ``eta'' is justified by their result according to which any 
first-order invertible self-adjoint differential operator $D$ can be canonically 
`suspended' to an invertible parametric differential operator 
$\cD \in \CL^1 (M,E; \R^{2k+1})$, whose higher eta invariant 
$\eta_{2k+1} (\cD)$ coincides with the spectral $\eta$-invariant $\eta (D)$. 
On the negative side, in contrast with Melrose's $\eta$-homomorphism,
these higher eta invariants are no longer additive on the multiplicative group of
invertible elements. Nevertheless,
the `defect of additivity' is purely symbolic and hence local.

The starting point for the developments that make the object of the present
paper was the fundamental observation that the higher $\eta$-invariants 
$\eta_{2k+1}$, when assembled together with symbolic corrections into 
\textit{higher divisor flows} $\DF_{2k+1}$, can be understood as the expression
of the Connes \textit{pairing} between the topological $K$-theory of the pair 
$\left(\CL^0(M,E; \R^{2k+1}), \CL^{-\infty} (M,E; \R^{2k+1})\right)$
and a certain canonical \textit{relative cyclic cocycle}, determined by the 
regularized \textit{graded trace} together with its symbolic coboundary. 
The first such invariant, for $k=0$, recovers Melrose's divisor flow,
whose essential properties such as homotopy invariance, additivity and 
integrality, thus acquire a conceptual explanation. Of course, the same 
properties are shared by the higher divisor flows  $\DF_{2k+1}$. Furthermore, 
this framework allows us to find the appropriate
even dimensional counterparts $\eta_{2k}$, $k > 0$ as well as their corresponding
higher divisor flows $\DF_{2k}$. Moreover, Theorem~\ref{bottcomp} confers them a
clear topological meaning, by showing that when
taken collectively the higher divisor 
flows $\DF_\bullet$ implement the natural Bott isomorphisms between the 
topological $K_\bullet$-groups of the pair 
$\left(\CL^0(M,E; \R^\bullet), \CL^{-\infty} (M,E; \R^\bullet)\right)$ and $\Z$, 
in a manner compatible with the suspension isomorphisms in both $K$-theory and 
in cyclic cohomology. 
 
As their basic properties indicate, and the very name given by  Melrose 
is meant to suggest, the divisor flows are close relatives of the spectral
flow. Our second main result gives a precise mathematical expression for this 
relationship, showing that spectral flow can be expressed as 
divisor flow via a natural Clifford
algebraic `suspension'. Because the standard complex Clifford 
representation comes in two `flavors', which
at the $K$-theoretical level provide the distinction between even and odd,
the result is formulated accordingly, in two separate
statements, Theorems \ref{sfodd} and \ref{sfeven}.

\tableofcontents

\section{Relative pairing in cyclic cohomology}
\subsection{Relative cyclic cohomology}
\label{Sec:RelCycCohom}
To establish the notation, we start by recalling in this section the 
definition of the relative homology and cohomology groups in terms
of pairs. We then specialize this description to the case of cyclic cohomology.

Consider a short exact sequence of chain complexes (over the field $\C$)
\begin{equation}
\label{Eq:ShExSeq}
  0 \longrightarrow (K_\bullet ,\partial_K ) 
  \overset{\kappa}{\longrightarrow} (A_\bullet ,\partial_A ) 
  \overset{\alpha}{\longrightarrow} (B_\bullet ,\partial_B ) 
  \longrightarrow 0,
\end{equation}
where the differentials $\partial_A, \partial_B, \partial_K$
are of degree $-1$. Put
\begin{equation}
  \widetilde K_k : = \operatorname{Cone} (\alpha)_{k+1}
  := A_k \oplus B_{k+1}, \quad
  \widetilde \partial :=  
  \left(  
  \begin{array}{cc}
     \partial_A & 0 \\
      - \alpha & - \partial_B 
  \end{array}
  \right)
\end{equation}
and
\begin{equation}
  \widehat B_k : = \operatorname{Cone} (\kappa)_{k+1} := 
  A_k \oplus K_{k-1}, \quad
  \widehat \partial :=  
  \left(  
  \begin{array}{cc}
     \partial_A & -  \kappa \\
     0 & - \partial_K 
  \end{array}
  \right).
\end{equation}
In other words, $\widetilde K_\bullet$ is the mapping cone of  $\alpha$ 
shifted by one degree, and $\widehat B_\bullet$ the mapping cone of $\kappa$.

Then $( \widetilde K_\bullet , \widetilde \partial )$ and
$( \widehat B_\bullet , \widehat \partial )$ 
are chain complexes, and one finds the following natural chain maps:
\begin{displaymath}
\begin{array}{lll}
  K_\bullet \rightarrow \widetilde K_\bullet , \:
  c_k \mapsto  (\kappa (c_k ),0), &
  \widehat B_\bullet \rightarrow B_\bullet, \:
  (a_k ,b_{k-1}) \mapsto \alpha (a_k), \\
  B_\bullet \overset{\iota}{\rightarrow} \widetilde K_{\bullet - 1}, \:
  b_k \mapsto (-1)^k (0,b_k), &
  A_\bullet  \overset{\iota}{\rightarrow} \widehat B_\bullet, \:
  a_k \mapsto (a_k,0), \\ 
  \widetilde K_\bullet \overset{\pi}{\rightarrow} A_\bullet, \:
  (a_k,b_{k+1}) \mapsto a_k ,&
  \widehat B_\bullet \overset{\pi}{\rightarrow} K_{\bullet -1}, \:
  (a_k,c_{k-1}) \mapsto (-1)^k c_{k-1}, 
\end{array}
\end{displaymath}
where we always assume that $a_l \in A_l$, $b_l \in B_l$ and $c_l\in K_l$.
Moreover, $K_\bullet \rightarrow \widetilde K_\bullet$ and
$\widehat B_\bullet \rightarrow B_\bullet$ are quasi-isomorphisms, and the
connecting morphism in the long exact homology sequence for 
\eqref{Eq:ShExSeq} becomes quite explicit. 
Namely, the long exact sequence reads
\begin{displaymath}
  \longrightarrow  H_k (\widetilde K_\bullet) \overset{\pi_*}{\longrightarrow}
  H_k(A_\bullet)  \overset{\alpha_*}{\longrightarrow}   
  H_k ( B_\bullet)  \overset{\iota_*}{\longrightarrow} 
  H_{k-1} (\widetilde K_\bullet) \longrightarrow  \, ,
\end{displaymath}
resp.
\begin{displaymath}
 \longrightarrow  H_k (K_\bullet)  \overset{\kappa_*}{\longrightarrow}  
  H_k(A_\bullet)  \overset{\iota_*}{\longrightarrow}   
  H_k ( \widehat B_\bullet)  \overset{\pi_*}{\longrightarrow} 
  H_{k-1} (K_\bullet) \longrightarrow  \, .
\end{displaymath}

We shall say in this case that the complex $\widetilde K_\bullet$
(resp.~$\widehat B_\bullet$)
encodes the relative homology of $A_\bullet \rightarrow B_\bullet$
(resp.~of $K_\bullet \rightarrow A_\bullet$).

In the dual situation one starts with a short exact sequence of cochain 
complexes
\begin{equation}\label{dualcochain}
  0 \longrightarrow (F^\bullet , d_F ) 
  \overset{\varepsilon}{\longrightarrow} (E^\bullet ,d_E ) 
  \overset{\delta}{\longrightarrow} (Q^\bullet ,d_Q ) 
  \longrightarrow 0,
\end{equation}
where the differentials $d_F, d_E, d_Q$ now have degree $+1$.
Then the cochain complex $(\widetilde Q^\bullet, \widetilde d)$, where
\begin{equation}\label{ML-G1.6}
  \widetilde Q^k : = E^k \oplus F^{k+1}, \quad
  \widetilde d :=  
  \left(  
  \begin{array}{cc}
     d_E & - \varepsilon \\
     0 & - d_F 
  \end{array}
  \right) 
\end{equation}
is quasi-isomorphic to  $(Q^\bullet,d_Q)$, and 
$(\widehat F^\bullet, \widehat d)$ with
\begin{equation}
  \widehat F^k : = E^k \oplus Q^{k-1}, \quad
  \widehat d :=  
  \left(  
  \begin{array}{cc}
     d_E & 0 \\
     - \delta & - d_Q 
  \end{array}
  \right), 
\end{equation}
is quasi-isomorphic to $(F^\bullet,d_F)$.
As above one obtains long exact cohomology sequences 
\begin{displaymath}
  \longrightarrow  H^k (F^\bullet) \overset{\varepsilon_*}{\longrightarrow}
  H^k(E^\bullet)  \overset{\iota_*}{\longrightarrow}   
  H^k (\widetilde Q^\bullet)  \overset{\pi_*}{\longrightarrow} 
  H^{k+1} (F^\bullet) \longrightarrow  \, ,
\end{displaymath}
resp.
\begin{displaymath}
 \longrightarrow  H^k (\widehat F^\bullet)  \overset{\pi_*}{\longrightarrow}  
  H^k(E^\bullet)  \overset{\delta_*}{\longrightarrow}   
  H^k (Q^\bullet)  \overset{\iota_*}{\longrightarrow} 
  H^{k+1} (\widehat F^\bullet) \longrightarrow  \, .
\end{displaymath}
Analogously to the homology case, we then say that the complex 
$\widetilde Q^\bullet$ (resp.~$\widehat F^\bullet$)
describes the relative cohomology of $F^\bullet \rightarrow E^\bullet$
(resp.~of $E^\bullet \rightarrow Q^\bullet$).

Recall that a cochain complex $E^\bullet$ is said to be {\it dual} 
to the chain complex $A_\bullet$, if there is a non-degenerate bilinear
pairing $\langle -, -\rangle : E^\bullet \times A_\bullet \rightarrow \C$ 
such that with respect to the pairing the differential $d_E$ is adjoint to
$\partial_A$. 
The main property of the complexes $\widetilde K_\bullet$ and 
$\widetilde Q^\bullet$ (resp.~$\widehat B_\bullet$ and $\widehat F^\bullet$) 
now is that they behave nicely under duality pairings. 
This is expressed in 
the following proposition the proof of which is straightforward.
\begin{proposition}
\label{Prop:DualPair}
Assume that in the short exact sequences \eqref{Eq:ShExSeq} and 
\eqref{dualcochain} the cochain complex $E^\bullet$ is dual to $A_\bullet$
and that $F^\bullet$ is dual to $B_\bullet$. Moreover, assume that
$\varepsilon$ is adjoint to $\alpha$. Then with respect to the pairing
\begin{equation}
\label{Eq:DualPair}
\begin{split}
  \langle - , -\rangle : \; &
  \widetilde Q^\bullet \times \widetilde K_\bullet \rightarrow 
  \C, \\
  & ( (\varphi_k,\psi_{k+1}), (a_k,b_{k+1})) \mapsto 
    \langle \varphi_k , a_k \rangle + \langle \psi_{k+1} , b_{k+1} \rangle 
\end{split}
\end{equation}
the complex $\widetilde Q^\bullet$ is dual to $\widetilde K_\bullet$, and this
pairing induces a bilinear pairing  
$H^k (\widetilde Q^\bullet) \times H_k (\widetilde K_\bullet) \rightarrow \C$.
Likewise, if $E^\bullet$ is dual to $A_\bullet$, $Q^\bullet$ dual to
$K_\bullet$, and $\delta$ adjoint to $\kappa$, 
then $\widehat F^\bullet$ is dual to $\widehat B_\bullet$.
\end{proposition}

We now specialize the above notions to give a description
of the relative cyclic homology and cohomology groups in the
`mapping cone setup'.

Recall that every unital $\C$-algebra $\cA$ (possibly endowed with a
locally convex topology) gives rise to a mixed 
complex $( C^\bullet (\cA),b,B)$, where 
$C^\bullet (\cA)$ is the Hochschild cochain complex, 
$b$ the Hochschild coboundary and $B$ Connes' boundary
(cf.~\textsc{Connes} \cite{Con:NDG}, \textsc{Loday} \cite{Lod:CH}). 
By definition, 
$C^k (\cA) = 
\big( \cA \otimes {\cA}^{\otimes k} \big)^* $,
and the operators $b$ and $B$ act on $\phi \in C^k (\cA)$ by 
\begin{equation}
\label{Eq:DefHb}
\begin{split}
  b \phi  (a_0 , \ldots , a_{k+1}) = \, 
  \sum_{j=0}^k &(-1)^j  \, \phi  ( a_0 , \ldots , a_j a_{j+1} , 
  \ldots , a_{k+1} ) \\
  & + (-1)^{k+1} \phi ( a_{k+1} a_0 , a_1 , \ldots , a_k ) ,
\end{split}   
\end{equation}
respectively
\begin{equation}
\label{Eq:DefCB}
\begin{split}
  B \phi& (a_0 , \ldots ,  a_{k-1}) \, = \, 
  \sum_{j=0}^{k-1} (-1)^{(k-1)j} \, \phi ( 1 , a_j , \ldots , a_{k-1}, 
  a_0 , \ldots ,  a_{j-1} )
 \\
  & - \,
  \sum_{j=0}^{k-1} (-1)^{(k-1)j} \, \phi (a_j , 1 , a_{j+1} , \ldots , a_k , 
  a_0 , \ldots , a_{j-1} ) .
\end{split}
\end{equation}

Consider now the double complexes 
$\mathcal BC^{\bullet,\bullet} (\cA)$ and 
$\mathcal BC_\text{\tiny per}^{\bullet,\bullet} (\cA)$.
They consist of the non-vanishing components 
$\mathcal BC^{p,q}(\cA) = C^{q-p} (\cA)$ for $q\geq p \geq 0$ 
resp.~$\mathcal BC_\text{\tiny per}^{p,q} (\cA) = C^{q-p} 
(\cA)$ 
for $q\geq p$, and have $B$ as horizontal and $b$ as vertical differential. 
The \emph{cyclic} resp.~\emph{periodic cyclic cohomology} groups of $\cA$ are then given
as follows:
\begin{displaymath}
  HC^\bullet (\cA) = 
  H^\bullet (\tot_{\bigoplus}^\bullet \mathcal BC^{\bullet,\bullet} 
  (\cA)) \: \text{ and } \:
  HP^\bullet (\cA) = H^\bullet (\tot_{\bigoplus}^\bullet 
  \mathcal BC_\text{\tiny per}^{\bullet,\bullet}(\cA)) .
\end{displaymath}
In both cases the differential on the total complex is $b+B$.
A (continuous) surjective homomorphism of algebras 
$\sigma : \cA \rightarrow \cB$ now induces a morphism of
mixed complexes
$\sigma^* : C^\bullet (\cB) \rightarrow C^\bullet (\cA) $.
Thus, we are in a `relative situation' as in \eqref{ML-G1.6}
and can express the corresponding
\emph{relative cyclic} resp.~\emph{periodic cyclic cohomology} accordingly. 
More precisely,  the relative cyclic cohomology coincides with the cohomology
of the total complex
\begin{displaymath}
  \big( \tot_{\oplus}^\bullet 
  \mathcal BC^{\bullet,\bullet} 
  (\cA) \oplus 
  \tot_{\oplus}^{\bullet+1} \mathcal BC^{\bullet,\bullet} 
  (\cB) , \widetilde{b+B}\,\big),
\end{displaymath}
where the differential is given by 
\begin{displaymath}
  \widetilde{b+B} = \left(
  \begin{array}{cc}
     b+B & -\sigma^* \\
     0 & -(b+B) 
  \end{array}
  \right).
\end{displaymath}
Explicitly, $\tot_{\oplus}^k  \, \mathcal B C^{\bullet,\bullet}
  (\cA) \oplus 
  \tot_{\oplus}^{k+1} \mathcal B C^{\bullet,\bullet}
  (\cB) \cong $
\begin{displaymath}
\begin{split}
  \cong
  \bigoplus_{p+q=k}  \mathcal B C^{p,q} (\cA) 
  \oplus  \mathcal B C^{p,q+1} (\cB)
  = \tot_\oplus^k \mathcal B C^{\bullet,\bullet}
  (\cA,\cB),
\end{split}
\end{displaymath}
where 
$\mathcal B C^{\bullet,\bullet} (\cA,\cB)$
is the double complex associated to the
relative mixed complex 
$(C^\bullet (\cA,\cB),\widetilde b,\widetilde B)$, which is 
given by 
$C^k (\cA,\cB) = C^k (\cA) \oplus C^{k+1} (\cB)$,
\begin{displaymath}
  \widetilde b =  \left(
  \begin{array}{cc}
     b & -\sigma^* \\
     0 & -b  
  \end{array}
  \right), \: \text{ and } \: 
  \widetilde B = \left(
  \begin{array}{cc}
     B & 0 \\
     0 & -B 
  \end{array}
  \right) .
\end{displaymath}
Hence the relative cyclic cohomology 
$HC^\bullet (\cA,\cB)$, resp.~the relative periodic 
cyclic cohomology 
$HP^\bullet (\cA,\cB)$ can be identified canonically with the 
cohomology of
$$
  \big( \tot_{\oplus}^\bullet \mathcal B C^{\bullet,\bullet}
  (\cA,\cB) , \widetilde b + \widetilde B\big)\quad
\text{resp.}  \quad
  \big( \tot_{\oplus}^\bullet \mathcal B C_\text{\tiny per}^{\bullet,\bullet}
  (\cA,\cB) , \widetilde b + \widetilde B\big).
$$
Of course, the ``staircase trick'' also works for the relative cyclic complex.
As a consequence, each class in $HC^k (\cA,\cB)$ has a 
representative 
$(\varphi,\psi)\in C^k_\lambda (\cA)\oplus C^{k+1}_\lambda (\cB)$
with $b\varphi = \sigma^* \psi$, where $C^\bullet_\lambda$ stands for
the subcomplex of cyclic cochains~\cite{Con:NDG}.

The preceding considerations can be dualized in an obvious fashion. Thus, 
$HC_\bullet (\cA,\cB)$ is the homology of 
$\big( \tot^{\oplus}_\bullet \mathcal B C_{\bullet,\bullet}
 (\cA,\cB) , \widetilde b + \widetilde B\big)$, 
where 
$\mathcal B C_{p,q} (\cA,\cB) = 
 \mathcal BC_{p,q} (\cA) \oplus 
 \mathcal BC_{p,q+1} (\cB)$, 
\begin{displaymath}
  \widetilde b =  \left(
  \begin{array}{cc}
     b & 0 \\
     -\sigma_* & -b  
  \end{array}
  \right), \: \text{ and } \: 
  \widetilde B = \left(
  \begin{array}{cc}
     B & 0 \\
     0 & -B 
  \end{array}
  \right).
\end{displaymath}
Likewise, the \emph{periodic cyclic homology} $HP_\bullet (\cA,\cB)$
is the homology of 
$\big( \tot^{\prod}_\bullet \mathcal B C^\text{\tiny per}_{\bullet,\bullet}
(\cA,\cB) , \widetilde b + \widetilde B\big)$, where 
$\mathcal B C^\text{\tiny per}_{p,q} (\cA,\cB) = 
 \mathcal BC^\text{\tiny per}_{p,q} (\cA) \oplus 
 \mathcal BC^\text{\tiny per}_{p,q+1} (\cB)$ and where
$ \widetilde b $, $\widetilde B$ are as above.

\subsection{Relative Chern character and relative pairing}
By Proposition \ref{Prop:DualPair}, the relative cyclic 
(co)homology groups inherit a natural pairing
\begin{equation}
\label{Eq:RelCycPair}
  \langle - , - \rangle_\bullet : \:
  HC_\bullet (\cA,\cB) \times HC^\bullet (\cA,\cB)
  \rightarrow \C ,
\end{equation}
which will be called the {\it relative cyclic pairing}. 

In Section 
\ref{Sec:DivFlowRelPair} below we shall express the divisor flow defined by
\textsc{Melrose} \cite{Mel:EIF} in terms of such a
relative pairing. In preparation for that, we recall below
the definition of the (odd) Chern character in periodic cyclic homology and 
the affiliated transgression formulas. 

Given a (Fr{\'e}chet) algebra $\cA$, and an element 
\[g\in \GL_\infty (\cA):=\lim\limits_{N\to\infty} \GL_N(\cA),\]
the odd Chern character is the following 
normalized periodic cyclic cycle: 
\begin{equation}
  \ch_\bullet (g) = \sum_{k=0}^\infty\, (-1)^k  k! 
  \tr_{2k+1} \big( (g^{-1} \otimes g)^{\otimes (k+1)} \big) , 
\end{equation}  
where $(g^{-1} \otimes g)^{\otimes j}$ is an abbreviation for the $j$-fold
tensor product \[ (g^{-1}\otimes g)\otimes\ldots\otimes (g^{-1}\otimes g),\]
and $\tr_k$ denotes the generalized trace map
$ \mathfrak{M}_N (\cA)^{\otimes k+1} \rightarrow 
  \cA^{\otimes k+1}$ (cf.~\cite[Def.~1.2.1]{Lod:CH}).

There is a transgression
formula (cf.~\textsc{Getzler}~\cite[Prop.~3.3]{Get:OCC})
for the odd Chern character of a smooth family of 
invertible matrices 
$g_s \in \GL_\infty (\cA)$, $s\in [0,1]$:
\begin{equation}
\label{Eq:ChTransgress}
  \frac{d}{ds} \ch_\bullet (g_s) = (b+B) \, 
  \slch_\bullet (g_s, \dot g_s),
\end{equation}
where the secondary Chern character $\slch_\bullet$ is defined by  
\begin{align}
  \slch_\bullet & (g , h ) = \tr_0 (g^{-1}h) + \\
    + & \sum_{k=0}^\infty (-1)^{k+1} k! 
    \sum_{j=0}^k \tr_{2k+2} \big( (g^{-1} \otimes g)^{\otimes (j+1)}\otimes
   g^{-1} h \otimes (g^{-1} \otimes g)^{ \otimes (k-j)} \big).\nonumber
\end{align}

Note that our sign convention differs from the one in 
\textsc{Getzler} \cite{Get:OCC}, 
since in Eqs.~\eqref{Eq:DefHb} and \eqref{Eq:DefCB} we have used the sign 
convention of e.g.~\cite{Lod:CH} for the definition of the 
operators $b$ and $B$.

The secondary Chern character fulfills the secondary transgression formula
Eq.~\eqref{Eq:SlChTransgress} for smooth two-parameter families of invertibles 
$g_{s,t}\in \GL_\infty (\cA)$, $s,t\in [0,1]$. 
Its proof follows by straightforward computation.  
\begin{equation}
\label{Eq:SlChTransgress}
  \frac{\partial}{\partial s} \slch_\bullet (g , \partial_t g ) 
  - \frac{\partial}{\partial t} \slch_\bullet (g , \partial_s g ) 
  = (b+B) \,  \sllch_\bullet (g, \partial_s g , \partial_t g),
\end{equation}
where
\begin{align}
  \sllch_\bullet & (g , h_1 , h_2 ) = - \tr_1 (g^{-1}h_1 \otimes g^{-1} h_2) -
    \\
    - & \sum_{k=0}^\infty (-1)^k k! 
    \sum_{j_1+j_2+j_3=k} \tr_{2k+3} \big( (g^{-1} \otimes g)^{\otimes (j_1+1)}
    \otimes g^{-1} h_1 \otimes \nonumber
    \\ &\hspace{40mm}
    \otimes (g^{-1} \otimes g)^{\otimes j_2} \otimes
    g^{-1} h_2 \otimes (g^{-1} \otimes g)^{ \otimes j_3} \big)  + \nonumber\\
    + & \sum_{k=0}^\infty (-1)^k k! 
    \sum_{j_1+j_2+j_3=k} \tr_{2k+3} \big( (g^{-1} \otimes g)^{\otimes (j_1+1)}
    \otimes g^{-1} h_2 \otimes\nonumber
    \\ &\hspace{40mm}
    \otimes (g^{-1} \otimes g)^{\otimes j_2} \otimes
    g^{-1} h_1 \otimes (g^{-1} \otimes g)^{ \otimes j_3} \big).\nonumber
\end{align}

The transgression formula Eq.~\eqref{Eq:ChTransgress} gives rise to a relative 
cyclic homology class as follows. As above, let 
$\sigma: \cA \rightarrow \cB$ be a continuous homomorphism of 
two locally convex topological algebras. 
Let us call an element $a\in\mathfrak{M}_N (\cA)$
{\it elliptic}, if $\sigma(a)$ is invertible, or in other 
words lies in $\GL_N(\cB)$. Moreover, given a family 
$(a_s)_{0\leq s \leq 1}$ of elements of $\mathfrak{M}_N (\cA)$, 
we say that $(a_s)_{0\leq s \leq 1}$ is an \emph{admissible elliptic path}, if 
each $a_s$ is elliptic, and $a_0$ and $a_1$ are both invertible.  

\begin{proposition}
\label{Prop:TransRelCyc}
 Let $(a_s)_{0\leq s \leq 1}$ be a smooth admissible elliptic path in 
 $\mathfrak{M}_N (\cA)$. Then the expression
 \begin{displaymath}
   \ch_\bullet \big( (a_s)_{0\leq s \leq 1} \big):=
   \Big( \ch_\bullet (a_1) - \ch_\bullet (a_0) , - \int_0^1 
   \slch_\bullet \big(\sigma(a_s) , \sigma(\dot a_s )\big) 
   \, ds \Big)
 \end{displaymath}
 is well-defined and defines a relative cyclic cycle.  

We then define the \emph{Chern character} of the admissible 
elliptic family $(a_s)_{0\leq s \leq 1}$ as the class of the relative  
cyclic cycle  $\ch_\bullet \big( (a_s)_{0\leq s \leq 1} \big)$.
\end{proposition}
\begin{proof}
  The transgression formula Eq.~\eqref{Eq:ChTransgress} entails
  \begin{displaymath}
     \sigma_* \big(\ch_\bullet (a_1) - \ch_\bullet (a_0) \big)
     = (b+B) 
     \int_0^1 \slch_\bullet 
     \big(\sigma (a_s) , \sigma (\dot a_s )\big) \, ds .
  \end{displaymath}
  Moreover, the odd Chern character is a cyclic cycle in $\cA$, and
  therefore  
  \begin{displaymath}
    (b+B) \big(\ch_\bullet (a_1) - \ch_\bullet (a_0) \big) =  0 , 
  \end{displaymath}
  which completes the proof.
\end{proof}

\subsection{Elliptic paths and relative $K$-theory}
  With the application alluded to above in mind, it is appropriate to
  elaborate at this point on invertibility and ellipticity 
  in topological algebras as well as on a representation of relative topological 
  $K$-theory by homotopy classes of elliptic paths.

  Assume that  $\cA$ is a topological algebra, 
  which in particular means that the product in $\cA$ is separately 
  continuous. Following the presentation in \textsc{Schweitzer} \cite{Schwe:SPLLF},
  we call  $\cA$ a {\it good} topological algebra, if 
  $\GL_1 (\widetilde{\cA}) \subset \widetilde{\cA}$ is open  and 
  the inversion is continuous, where $\widetilde{\cA}$ denotes
  the smallest unital algebra containing $\cA$.
  By \textsc{Swan} \cite[Lem.~2.1]{Swa:TEPM},  $\cA$ being a good 
  topological algebra implies that $\mathfrak{M}_N (\cA)$ is good as well for 
  every $N \in \N^*$. Moreover, if $\sigma : \cA \rightarrow \cB$ 
  is a continuous homomorphism of topological algebras with  $\cA$ a good 
  topological algebra, then $\cB$ is good as well by \textsc{Gramsch} 
  \cite[Bem.~5.4]{Gra:RISOA}.  Thus, in this situation, 
  the space $\Ell_N (\widetilde{\cA})$ of elliptic 
  $N\times N$ matrices with entries in $\widetilde{\cA}$ is an open 
  subset of $\mathfrak{M}_N (\widetilde{\cA})$. If $\cA$ is additionally
  locally convex, this implies in particular that for every continuous family 
  $(a_s)_{0\leq s \leq 1}$ of elliptic elements in 
  $\mathfrak{M}_N (\widetilde{\cA})$, one can find a smooth path
  $(b_s)_{0\leq s \leq 1}$ which, in $\Ell_N (\widetilde{\cA})$, 
  is homotopic to  $(a_s)_{0\leq s \leq 1}$ relative the endpoints.  
  Moreover, using \cite[Bem.~5.4]{Gra:RISOA} again, one can even show 
  that there exists a continuous family $(q_s)_{0\leq s \leq 1}$ of 
  {\it parametrices} for 
  $(a_s)_{0\leq s \leq 1}$, which means that both $a_s q_s - I$ and 
  $q_s a_s -I$ lie in the kernel of $\sigma$ for each $s\in [0,1]$.  

  Recall that a locally convex topological algebra is called 
  {\it m-convex}, if its topology is generated by submultiplicative seminorms. 
  Next, a Fr\'echet subalgebra 
  $ \iota : \cA \hookrightarrow  \cA'$ of some m-convex 
  Fr\'echet algebra $\cA'$ is called  
  {\it closed under holomorphic functional calculus}, if  for every 
  $A \in \widetilde{\cA}$ and every complex function $f$ holomorphic 
  on a neighborhood of the  $\widetilde{\cA'}$-spectrum of $A$, 
  the element $f(A)\in \widetilde{\cA'}$ lies in the algebra
  $\widetilde{\cA}$. If  $\cA'$ is a Banach algebra 
  (resp.~$C^*$-algebra) and $\cA$ is dense in $\cA'$, then one  
  often calls $\cA$ a {\it local Banach algebra} 
  (resp.~{\it local $C^*$-algebra}) or briefly {\it local} in $\cA'$, 
  and denotes $\cA'$ as $\bar{\cA}$, the {\it completion} of 
  $\cA$.   
  By a result going back to the work of \textsc{Gramsch} \cite{Gra:RISOA} and 
  \textsc{Schweitzer} \cite{Schwe:SPLLF}, the matrix algebra $\mathfrak{M}_N (\cA)$ 
  over a Fr\'echet algebra $\cA$ is local, if and only if  $\cA$ is. 
  This implies in particular the following stability result.
\begin{proposition}
\label{Thm:StabKThe}
\textup{(Cf.~\cite[VI.~3]{Con:ATICPAA})}
  If $\cA $ is a unital Fr\'echet algebra which is local in a Banach 
  algebra $\bar{\cA}$, then 
  $ \iota: \cA \hookrightarrow \bar{\cA}$ induces an isomorphism 
  of topological $K$-theories:
  \begin{displaymath}
    \iota_* : K_\bullet (\cA) \cong K_\bullet (\bar{\cA}) .
  \end{displaymath}
\end{proposition}

Let us now briefly recall the construction of the relative $K_1$-group associated 
to a continuous unital homomorphism $\sigma : \cA \rightarrow \cB$ 
of unital Fr\'echet algebras. The group $K_1 (\cA, \cB)$ is defined 
as the kernel of the map 
${\operatorname{pr}_2}_* :K_1 (D) \rightarrow  K_1 (\cA)$
induced by the projection 
$\operatorname{pr}_2 : D \rightarrow \cB$ of the double 
\[ 
  D = \big\{ (a,b ) \in \cA \times \cA \mid a-b \in 
  \cJ:= \Ker \sigma 
  \big\} 
\] 
onto the second coordinate. Since 
$K_1 (\cA )$ is the quotient of  $\GL_\infty  (\cA)$
by the connected component of the identity, it is clear that   
$K_1 (\cA, \cB)$ can be naturally identified with the
quotient of
\begin{displaymath}
  G := \big\{ a \in \GL_\infty (\cA ) \mid a  = I + S
  \text{ with $ S \in \cJ $ } \big\}
\end{displaymath}
by $G_0$, the connected component of the identity in $G$. 
By excision in (topological) $K$-theory, one knows that 
$K_1 (\cA,\cB)$ is naturally isomorphic to $K_1 (\cJ)$. 
We will present another identification of $K_1 (\cA,\cB)$, 
namely with a certain fundamental group of the space 
\begin{displaymath}
   \Ell_\infty (\cA) := \lim_{N\rightarrow \infty} \Ell_N (\cA) 
\end{displaymath} 
of elliptic elements of $\cA$. Note that as an inductive limit this 
space carries a natural topology inherited from the Fr\'echet topology on 
$\cA$.

Before stating the result we need to fix some more notation. For a topological 
space $X$ and a subspace $Y\subset X$ we denote by $\Omega(X,Y)$ the set of 
continuous paths in $X$ with endpoints in $Y$. If $y_0\in Y$ is a basepoint 
in $Y$ we denote by $\pi_1(X,Y;y_0)$ the relative homotopy set of homotopy 
classes of paths with initial point $y_0$ and endpoint in $Y$. 
By $\pi_1(X,Y)$ we denote the homotopy classes of paths $\sigma:[0,1]\to X$
with $\sigma(0),\sigma(1)\in Y$. Note that for a homotopy of paths the 
endpoints may vary in $Y$. $\pi_1(X,Y)$ together with concatenation of paths 
is naturally a groupoid. For definiteness we denote the concatenation 
of paths by $*$.

On $\pi_1(\Ell_\infty(\cA), \GL_\infty (\cA))$, the set of 
homotopy classes of admissible elliptic paths,  we even have a monoid 
structure, $\cdot$, induced by pointwise multiplication. 
In fact, the two products are not independent:
consider paths $f,g\in\Omega(\Ell_\infty(\cA),\GL_\infty(\cA))$ and put
\begin{equation}
  \varphi(s)=\begin{cases} 2s,&0\le s\le 1/2,\\ 1,&1/2\le s\le 1\end{cases}, 
  \quad
  \psi(s)=\begin{cases} 0,& 0\le s \le 1/2,\\ 2s-1,& 1/2\le s\le 1.\end{cases}
\end{equation}
Then $H(s,t)=f((1-t)s+t\varphi(s))g((1-t)s+t\psi(s))$
is a homotopy which shows that
\begin{equation}
  \label{MLEq:homotopyconcatenation}
  f\cdot g\simeq (f\cdot g(0))* (f(1)\cdot g).
\end{equation} 

Denote by $f_-(s):=f(1-s)$ the inverse of the path $f$ with respect to the 
concatenation product. 
Then $f*f_-$ is homotopic to the constant path $f(0)$. In view of 
\eqref{MLEq:homotopyconcatenation} we have 
$f\cdot(f(1)^{-1} f_-)\simeq f*f_-\simeq f(0)$
and thus we have shown that for every path $f$ there is a path $g$ such that
$f\cdot g$ is homotopic to a constant path.

Finally, we call a path $f:[0,1]\longrightarrow \GL_\infty(\cA)$ 
\emph{degenerate}. A path is degenerate, if and only if it is homotopic to a 
constant path. The homotopy classes of degenerate paths form a submonoid of 
$\pi_1(\Ell_\infty(\cA),\GL_\infty(\cA))$.
If we mod out this submonoid of $\pi_1(\Ell_\infty(\cA),\GL_\infty(\cA))$
we obtain a group which we denote by 
$\widetilde\pi_1(\Ell_\infty(\cA),\GL_\infty(\cA))$.
Observe that by Eq.~\eqref{MLEq:homotopyconcatenation} the products
$*$ and $\cdot$ coincide on 
$\widetilde\pi_1(\Ell_\infty(\cA),\GL_\infty(\cA))$.
Finally we note that thanks to the fact that we are working with paths in a 
stable algebra the pointwise product in  
$\pi_1(\Ell_\infty(\cA),\GL_\infty(\cA))$
may as well be represented by direct sums since $f\cdot g$ is homotopic to
\begin{equation}
  \begin{pmatrix}
      f&0\\ 0& g
    \end{pmatrix}.
\end{equation}

\begin{lemma}\label{ML-Lemma1.4}
The natural map
\begin{equation}\pi_1(\Ell_\infty(\cA),\GL_\infty(\cA);I)\longrightarrow
\widetilde\pi_1(\Ell_\infty(\cA),\GL_\infty(\cA))
\end{equation}
is an isomorphism.
\end{lemma}
\begin{remark} It is worth pointing out that
  this fact does not have an exact analogue in the even case (see 
  Section \plref{evencase} below). The reason is that $\GL_\infty(\cA)$ has 
  a distinguished base point while the space of idempotents in 
  $\mathfrak{M}_\infty(\cA)$ does not.
\end{remark}

\begin{proof} 
If $f\in\Omega(\Ell_\infty(\cA),\GL_\infty(\cA))$, then the class of the path 
$f$ in $\widetilde\pi_1(\Ell_\infty(\cA),\GL_\infty(\cA))$ is represented by 
the path $f(0)^{-1}f$ starting at $I$. This proves surjectivity. To prove 
injectivity consider a path
$f\in\Omega(\Ell_\infty(\cA),\GL_\infty(\cA))$ starting at $I$ such that $f$ 
is homotopic to a degenerate path. Then there is a map 
$H:[0,1]^2\to \Ell_\infty(\cA)$ such that $H(s,0)=f(s)$, $H(s,1)$ is a fixed 
element of $\GL_\infty(\cA)$, and $H(0,t), H(1,t)\in\GL_\infty(\cA)$.
Hence $\widetilde H(s,t):=H(0,t)^{-1}H(s,t)$ is a homotopy of paths starting 
at $I$ which implements a homotopy between $f$ and the constant $I$.
\end{proof}
\begin{theorem}
\label{Thm:RelKEll}
  Let  $\cA$ and $\cB$ be two unital Fr\'echet algebras and 
  \begin{displaymath}
    0 \longrightarrow \cJ \longrightarrow  \cA 
    \overset{\sigma}{\longrightarrow} \cB \longrightarrow 0
  \end{displaymath}
  an exact sequence of Fr\'echet algebras and unital homomorphisms such that   
  $\cA$ (and hence $\cB$) is a good Fr\'echet algebra. 
  Moreover, assume that $\cA$ and $\cB$ are both local Banach algebras 
  such that $\sigma$ extends to a continuous homomorphism 
  $\bar{\cA} \rightarrow \bar{\cB}$.
  Then  the following holds true:
  \begin{enumerate}
  \item[\textup{(1)}] 
     For each $a_0 \in \GL_\infty (\cA)$, the inclusion 
     $\cA\hookrightarrow \bar{\cA}$ induces a natural bijection
     \begin{displaymath}
        \pi_1 \big( \Ell_\infty (\cA) , \GL_\infty (\cA) ; a_0 \big) 
        \cong
        \pi_1 \big( \Ell_\infty (\bar{\cA}) , \GL_\infty (\bar{\cA}) ; 
        a_0 \big).
     \end{displaymath} 
  \item[\textup{(2)}] 
      The canonical homomorphism
      \begin{displaymath}
      \begin{split}
         \kappa: \: K_1 (\cA, \cB) \, & \rightarrow  
        \pi_1 \big( \Ell_\infty (\cA) , \GL_\infty (\cA) ; I \big), \\
        \text{induced by the assignment} \\
         I + S \, & \mapsto \big[ [0,1] \ni s \mapsto  I + s S \in  
         \Ell_\infty (\cA) \big]
      \end{split}
      \end{displaymath} 
      is an isomorphism, which therefore induces a canonical isomorphism
$K_1(\cA,\cB)\longrightarrow \widetilde\pi_1(\Ell_\infty(\cA),\GL_\infty(\cB))$.
  \end{enumerate}
\end{theorem}
\begin{proof}
The proof of the first claim will only be sketched, since it also follows from 
the second and stability of $K$-theory for local Banach algebras. 
The locality implies in particular that 
$\GL_N (\cA) = \GL_N (\bar{\cA}) \cap \mathfrak{M}_N (\cA)$ 
for each $N$ (see e.g.~\cite{Schwe:SPLLF}).
Likewise, one shows 
$\Ell_N (\cA) = \Ell_N (\bar{\cA}) \cap \mathfrak{M}_N (\cA)$. 
Since
$\Ell_N (\bar{\cA})$ is open in $\mathfrak{M}_N (\bar{\cA})$, and 
the topology of $\bar{\cA}$ is generated by a complete norm, every 
path in $\Ell_N (\bar{\cA})$ starting at $a_0$ is homotopic in  
$\Ell_N (\bar{\cA})$ to a path with values in $\Ell_N (\cA)$.  
This proves (1). 

The second claim is essentially a consequence of the fact that under 
the assumptions made
(notably locality of $\cA$ and $\cB$), there exists for each 
$b \in \GL_N (\cB)_I$, the connected component of the identity in 
$\GL_N (\cB)$, an invertible lift, i.e.~an element 
$a\in \GL_N (\cA)_I$ with $\sigma (a)=b$
(see for example \textsc{Blackadar} \cite[Cor.~3.4.4]{Bla:TOA} or 
\textsc{Gramsch} \cite[Bem.~5.4]{Gra:RISOA}). 
This lifting property holds for continuous paths as well, since with 
$\cA$ being local, 
$\cC^\infty ([0,1], \cA) \cong \cC^\infty ([0,1])\hat{\otimes} \cA$ 
is local as well, and likewise for $\cB$. Note also that by nuclearity
of $\cC^\infty ([0,1])$ the sequence 
\begin{displaymath}
  0 \longrightarrow \cC^\infty ([0,1], \cJ) \longrightarrow 
  \cC^\infty ([0,1], \cA) \longrightarrow 
  \cC^\infty ([0,1], \cB) \longrightarrow 0
\end{displaymath}
is exact as well, and that the closure of $\cC^\infty ([0,1], \cA)$ 
is $\cC ([0,1], \bar{\cA})$.

After these preliminary considerations let us now prove surjectivity of 
$\kappa$.
Let $(a_s)_{0\leq s \leq 1}$ be a smooth path in $\Ell_N (\cA)$ 
which starts at the identity and satisfies $a_1 \in \GL_N(\cA)$. 
By the preceding remarks there exists an invertible lift of 
$(\sigma (a_s))_{0\leq s \leq 1}$, that means a smooth path 
$(b_s)_{0\leq s \leq 1}$ of elements of $\GL_N (\cA)$ such that 
$\sigma(a_s) = \sigma (b_s) $ for all $s$. Hence the path 
$( a_s b_s^{-1} )_{0\leq s \leq 1}$
consists of elements of the form $I + S_s$, where $S_s \in \cJ$. 
Moreover, since the homotopy class of $(b_s)_{0\leq s \leq 1}$ is trivial,
$(a_s)_{0\leq s \leq 1}$ and  $( a_s b_s^{-1} )_{0\leq s \leq 1}$ are 
homotopic. Consider now the homotopy
\begin{displaymath} 
  h(s,t) = t \, a_s b_s^{-1} + (1-t) \big( I + s(a_1 b_1^{-1} -I ) \big ),
\end{displaymath}
where $s\in [0,1]$ and $t\in [0,1]$ is the homotopy parameter.
Then one has 
\begin{enumerate}\romlabel
  \item
     $h(0,t) =I$, 
  \item
    $ h(1,t) = a_1b_1^{-1}$, which is invertible, 
  \item 
    $h(s,1) = a_sb_s^{-1}$, and  
  \item
     $h(s,0) = I + s S_1 $ with $S_1 := (a_1b_1^{-1} -I) \in \cJ$ .   
 \end{enumerate}
 Hence $(a_s)_{0\leq s \leq 1}$ is homotopic to the path $h(.,0)$, which is 
equal to  $\kappa (I+S_1)$. This shows surjectivity.

 Let us finally show injectivity of $\kappa$. Assume that $\kappa (I+S)$ has 
 trivial homotopy class.
 This means that there exists a homotopy $h$ in $\Ell_N (\cA)$ relative 
 $\GL_N (\cA)$
 with $h(s,0) = I+sS$, $h(s,1) =I$ and $h(0,t) =I$. 
 Since $\sigma \circ h$ is a homotopy of invertible elements in 
 $\GL_N(\cB)$, there exists an invertible lift to $\cA$, 
 i.e.~a homotopy $g$ in $\GL_N (\cA)$ such that 
 $\sigma \circ g = \sigma \circ h$ and
 \begin{enumerate}\romlabel
  \item
    $g(0,t)=g(s,1) =I$ for all $s,t\in [0,1]$.
 \end{enumerate}
 Furthermore, one has 
 \begin{enumerate}\romlabel\setcounter{enumi}{1}
  \item
    $g(s,0)=I +T_s$ with $T_s\in J$.
  \end{enumerate}   
 Consider now the homotopy 
 $H(s,t):= h (s,t) \, g^{-1} (s,t) \, (I+ T_{(1-t)s})$.
 Since $   h (s,t) \, g^{-1} (s,t) -I \in \cJ$ for all $s,t$, one 
 concludes that
\begin{enumerate}\romlabel\setcounter{enumi}{2}
  \item
     $H(0,t) =I$, 
  \item
    $ H(s,0) = h(s,0) = I +sS$,  
  \item 
    $H(s,1) = I$, and
  \item 
    $H(s,t) -I \in \cJ$ for all $s,t$. 
 \end{enumerate} 
 Therefore, $H(s,.)$ gives rise to a homotopy in $K_1(\cA,\cB)$ 
 between $I+S$  and $I$. Hence $\kappa$ is injective, and the theorem follows.
\end{proof}

Using the preceding theorem we can now show that the Chern character from the 
previous section is defined even on 
$\pi_1 \big(\Ell_\infty (\cA), \GL_\infty (\cA)\big)$, 
the fundamental groupoid of the space of elliptic elements of the
algebra $\cA$ relative the invertible ones.
\begin{theorem}
\label{Thm:NatChern}
 Under the assumptions of the preceding theorem 
 let $(a_{s,t} )_{0\leq s,t \leq 1}$
 be a smooth family in $\mathfrak{M}_N(\cA)$ such that for each fixed $t$ 
 the family $(a_{s,t} )_{0\leq s \leq 1}$ is a smooth admissible path. 
 Then $ \ch_\bullet \big( (a_{s,1} )_{0\leq s \leq 1} \big)$ and
 $ \ch_\bullet \big( (a_{s,0} )_{0\leq s \leq 1} \big)$ are
 homologous relative cyclic cycles, and therefore 
  the Chern character $\ch_\bullet$ is homotopy
 invariant and descends to a map on 
 $\pi_1 \big( \Ell_\infty (\cA) , \GL_\infty (\cA)\big)$.
 Moreover, if $\cA$ and $\cB$ are both local Banach algebras, then
 $$\ch_\bullet \circ \kappa : K_1 (\cA, \cB) \rightarrow HP_1
 (\cA, \cB)
 $$ 
 coincides, via the canonical identification  
 $K_1 (\cJ) \cong K_1  (\cA, \cB)$, where 
 $\cJ = \Ker \sigma$, with the standard Chern character in cyclic homology.
\end{theorem}
\begin{proof}
The secondary transgression formula Eq.~\eqref{Eq:SlChTransgress} entails 
that 
\begin{displaymath}
\begin{split}
   \ch_\bullet & \big( (a_{s,1} )_{0\leq s \leq 1} \big)
   - \ch_\bullet \big( (a_{s,0} )_{0\leq s \leq 1} \big) \\
   & = \big( \widetilde b + \widetilde B \big) 
   \Big( \slch \big( a_{s,1}, \partial_t a_{s,1} \big) 
   - \slch \big( a_{s,0}, \partial_t a_{s,0} \big) , \\
   & \hspace{30mm}- \int_0^1 \sllch \big( \sigma (a_{s,t} ) ,\sigma 
   (\partial_s a_{s,t}), \sigma (\partial_t a_{s,t})\big) ds \Big),
\end{split}
\end{displaymath}
which proves the first part of the theorem. To show the second part let
$[I +S] \in K_1(\cJ)$ and check that
\begin{displaymath}
  \ch_\bullet \big( \kappa ([I+S]) \big) =
  \big( \ch_\bullet (I) -\ch_\bullet (I+S), 0) .
\end{displaymath}
This gives the claim.
\end{proof}

The higher homotopy groups 
$\pi_n \big( \Ell_\infty (\cA),\GL_\infty (\cA); I\big)$, 
$n\geq 1$, satisfy a natural Bott periodicity property.
In order to state it, let us assume to be given a homotopy
class 
$[\varrho] \in \pi_2 \big( \Ell_\infty (\cA) , \GL_\infty (\cA) ;I \big)$,
represented by a continuous map
\begin{displaymath}
  \varrho : [0,1]^2  \rightarrow \Ell_N (\cA)  ,
\end{displaymath}
with the additional property  that 
\begin{displaymath}
 \varrho \big( [0,1] \times \{ 0 \} \big) = I \quad \text{and} \quad 
 \varrho \big(\partial [0,1]^2 \big) \subset \GL_N (\cA).
\end{displaymath}
Obviously, $\varrho$ is homotopic to a continuous map 
$\widetilde \varrho : [0,1]^2  \rightarrow \Ell_N (\cA) $
satisfying
\begin{displaymath}
 \widetilde \varrho \big( [0,1] \times \{ 0 \}\cup \{0,1\}\times [0,1]\big) = I 
 \quad \text{and} \quad \widetilde\varrho\big(\partial [0,1]^2\big)\subset \GL_N (\cA).
\end{displaymath}
Since $\widetilde{S\cA}$, the suspension of $\cA$ with unit adjoined, naturally 
coincides with $\bigl\{f\in\cC^0 (S^1,\cA)\,\big|\, f(1)\in \C 1_{\cA}\bigr\}$, the path $\widetilde \varrho (\cdot , t)$ is
an element of $\widetilde{S\cA}$ for each $t\in [0,1]$.
Hence the map $[0,1] \ni t \mapsto \widetilde \varrho (\cdot , t)$ defines
a homotopy class in $\pi_1 \big( \Ell_\infty (\widetilde {S \cA}), 
\GL_\infty ( \widetilde{ S\cA} );I \big) $
which we denote by $\lambda ( \varrho)$.  Clearly, $\lambda$ factors
through $\pi_2 \big( \Ell_\infty (\cA) , \GL_\infty (\cA) ;I \big)$ and
one checks immediately that the resulting map 
\begin{equation}
  \lambda: \pi_2 \big( \Ell_\infty (\cA) , \GL_\infty (\cA) ;I \big)
  \rightarrow 
  \pi_1 \big( \Ell_\infty (\widetilde {S \cA}), 
  \GL_\infty ( \widetilde{ S\cA} );I \big)
\end{equation}
is an isomorphism. Using suspension in $K$-theory one then finds
\begin{equation}\begin{split} 
  \pi_2 \big( \Ell_\infty (\cA) , \GL_\infty (\cA) ;I \big)
   & \cong \pi_1 \big( \Ell_\infty (\widetilde {S \cA}) , 
   \GL_\infty ( \widetilde{ S\cA} );I \big) \cong \\
   &\cong  K_1 (S\cA, S\cB) \cong K_0 (\cA, \cB).
 \end{split}
\end{equation}
By iteration of this argument and using Bott periodicity in $K$-theory,
one obtains the corresponding periodicity property of the higher homotopy 
groups in the relative setting.

\begin{theorem}
 Let  $\sigma : \cA  \rightarrow \cB$ be a surjective
 homomorphism of good Fr\'echet algebras such that the assumptions of
 Thm.~\ref{Thm:RelKEll} are satisfied. Then  
 \begin{displaymath}
  \pi_{2n-i} \big( \Ell_\infty (\cA) , \GL_\infty (\cA);I \big) 
  \, = \,  K_i  (\cA, \cB) , \text{ for } i=0,1, \text{ and all } n\ge 1.
 \end{displaymath}
 \end{theorem}

\subsection{Relative cycles and their characters}
\label{sec1.3}
Connes' concept of a cycle over an algebra has a natural extension to
the relative case.  In what follows
   ${\cA}, {\cB}$ are unital 
(locally convex topological) $\C$-algebras and  
$\sigma:{\cA}\to {\cB}$ 
is a (continuous) unital and surjective homomorphism. We recall from
\cite[Part II.~\S 3]{Con:NDG} (cf. also~\cite[Sec.~2]{Gor:CCE}) the definition
of a `chain', which we prefer to call ``relative cycle'' here.

\begin{definition}
\label{defcyclewithboundary}
A {\it relative cycle} of degree $k$ over $(\cA,\cB)$ 
consists of the following data: 
\begin{enumerate}
\item 
  differential graded unital algebras $(\Omega , d)$ and 
  $(\partial \Omega , d)$ over $\cA$ resp.~$\cB$ together 
  with a  surjective unital homomorphism $r:\Omega\to \partial\Omega$ of 
  degree $0$,
\item
  unital homomorphisms $\varrho_{\cA}:\cA\to\Omega^0$ and
  $\varrho_{\cB}:\cB\to\partial\Omega^0$ such that
  $r\circ\varrho_{\cA}=\varrho_{\cB}\circ \sigma$,
\item 
 a graded trace $\int$ on $\Omega$ of degree $k$ such that
 \begin{equation}
       \int  d\omega=0 \, ,  \quad \text{whenever} \quad r( \omega) = 0.
\end{equation} 
\end{enumerate}
The graded trace $\int$ induces a unique closed graded trace $\int'$ on 
$\partial \Omega$ of degree $k-1$, such that Stokes' formula
\begin{equation}
  \int d\omega=\int' r\omega \quad \text{ for all $\omega\in\Omega$}
\end{equation}
is satisfied.

A relative cycle will often be denoted in this article as a tuple
$C = \big( \Omega,\partial \Omega,r,\rho_{\cA},\rho_{\cB}, 
\int,\int' \big)$ or, more compactly, as
$C= \big( \Omega,\partial \Omega,r,\int,\int' \big)$.  
The boundary $(\partial \Omega,d,\int')$ is just a cycle over the algebra 
${\cB}$. 
For $(\Omega,d,\int)$, this is in general not the case, unless the trace $\int$ 
is closed.   
\end{definition}

Note that, for the unit $1_\Omega$ of $\Omega$, the Leibniz rule implies
$d 1_\Omega=0$. By definition we consider only cycles with $\varrho_{{\cA}}$
unital. This simplifies our considerations because of 
$d\varrho(1_{\cA})=d 1_\Omega=0$.
(When dealing with non-unital $\varrho_{\cA}$ one has to accept
the disturbing fact that $d\varrho_{\cA}(1)$ might be nonzero.)
If no confusion can arise, we will omit from now on the subscript of $\varrho$.

We next define  the \emph{character} of a relative cycle $C$ by means of the
following proposition.

\begin{proposition}
\label{ML-prop3.1}\textup{(Cf.~\cite[Sec.~2]{Gor:CCE})} 
Let $C$ be a relative cycle  of degree $k$ over $(\cA,\cB)$. 
Define 
$(\varphi_k,\psi_{k-1})\in C^k({\cA})\oplus C^{k-1}({\cB})$ 
as follows:
\begin{align}
    &\varphi_k (a_0,\ldots,a_k):=\frac{1}{k!}\,
     \int\varrho(a_0)d\varrho(a_1)\ldots d\varrho(a_k),\\
    &\psi_{k-1}(b_0,\ldots,b_{k-1}):=  \frac{1}{(k-1)!}\,
     \int'\varrho(b_0)d\varrho(b_1) \ldots d\varrho(b_{k-1}).
\end{align}
Then $\character C:=(\varphi_k,\psi_{k-1})$ is a relative cyclic cocycle
in  $\tot_\oplus^k\mathcal B C^{\bullet,\bullet} 
   ({\cA},{\cB})$, called the character of the relative cycle $C$.
\end{proposition}
\begin{proof}
Since $(\partial \Omega,d,\int')$ is a cycle over $\cB$, we know
by the work of \textsc{Connes} \cite{Con:NDG} that $\psi_{k-1}$ is a cyclic
cocycle over $\cB$, which in particular means that $(b+B)\psi_{k-1}=0$.
Hence it remains to show that
\begin{equation}\label{ML-G1.14}
     (b+B)\varphi_k=\sigma^*\psi_{k-1}.
\end{equation}
Eq.~\eqref{ML-G1.14} in fact means
\begin{displaymath}
    b\varphi_k=0, \quad B\varphi_k=\sigma^*\psi_{k-1}.
\end{displaymath} 
For simplicity we omit $\varrho$ in the notation.
Now, $b$-closedness follows from the trace property. 
Namely, for $a_0,\ldots, a_{k+1}\in {\cA}$ applying Leibniz rule to 
$d(a_j a_{j+1})$ we find
\begin{equation}
  \begin{split}
    \sum_{j=1}^k &(-1)^j \int a_0 da_1\cdot\ldots\cdot d(a_j a_{j+1})
    \cdot\ldots\cdot da_{k+1} = \\
    &= - \int a_0 a_1 da_2\cdot\ldots\cdot da_{k+1}+(-1)^k 
    \int a_0 da_1\cdot\ldots\cdot da_k a_{k+1}
  \end{split}
\end{equation}
and hence by the trace property $b \varphi_k=0$. 
Using $d1=0$ and Stokes formula for $C$ we find 
\begin{align}
   B \varphi_k & (a_0,\ldots,a_{k-1})\nonumber\\
     = \, & \sum_{j=0}^{k-1} (-1)^{(k-1) j} \, \varphi_k 
       ( 1 , a_j , \ldots , a_{k-1} , a_0 , \ldots ,  a_{j-1} ) \\
     & \qquad -  \sum_{j=0}^{k-1} (-1)^{(k-1)j} \, \varphi_k 
     (a_j , 1 , a_{j+1} , \ldots , a_k , a_0 , \ldots , a_{j-1} ) \displaybreak[3] \nonumber\\ 
     = \, & \frac{1}{k!} \sum_{j=0}^{k-1} (-1)^{(k-1)j} \int 
     da_j\cdot\ldots\cdot da_{k-1} \, da_0 \cdot \ldots \cdot da_{j-1}\nonumber\\
     = \, & \frac{1}{(k-1)!} \int d \big( a_0 da_{1}\cdot\ldots\cdot 
     da_{k-1}\big) \nonumber\\
     = \, & \frac{1}{(k-1)!} \int' r(a_0 )dr(a_{1})\cdot\ldots\cdot dr (a_{k-1}) \nonumber\\
     = \, & \psi_{k-1}(\sigma (a_0) ,\ldots,\sigma (a_{k-1}) ),\nonumber
\end{align}
where we have used that $\int$ is a graded trace.
 \end{proof}

\subsection{Divisor flow associated to an odd relative cycle}
\label{Sec:AbsDivisor}
In the situation of Section \ref{sec1.3}, consider a relative cycle
$C$ of degree $2k+1$ over $(\cA,\cB)$. 
Denote by $(\varphi_{2k+1},\psi_{2k})$ the character of  $C$ as 
defined above, and recall that as a consequence of Proposition 
\ref{Prop:DualPair} we have a convenient representation of the 
natural pairing \eqref{Eq:RelCycPair} in relative cyclic (co)homology. 

\begin{definition} 
The ({\it odd}) {\it divisor flow} of the smooth admissible elliptic path 
$(a_s)_{0 \leq s\leq 1}$  with respect to the odd relative cycle $C$ is 
defined to  be the relative pairing between 
$\ch_\bullet \big( (a_s)_{0\leq s \leq 1} \big)$ and 
the character $\character C=(\varphi_{2k+1},\psi_{2k})$, i.e.
\begin{equation}\label{Eq:dfodd}
 \begin{split}
    \DF_C \big(  &(a_s)_{0\leq s\leq 1} \big) := 
    \DF \big( (a_s)_{0\leq s\leq 1} \big) := \\
    := \, & \frac{1}{(-2\pi i)^{k+1}}\langle (\varphi_{2k+1},\psi_{2k}),
    \ch_\bullet ((a_s)_{0\le s\le 1})\rangle \\
    = \, &\frac{1}{(-2\pi i)^{k+1}} \big( \langle 
    \varphi_{2k+1},\ch_\bullet(a_1)\rangle -
    \langle \varphi_{2k+1}, \ch_\bullet(a_0)\rangle\big) \\  
    & -\frac{1}{(-2\pi i)^{k+1}}\langle \psi_{2k},\int_0^1 
    \slch_\bullet(\sigma (a_s), \sigma (\dot a_s))ds \rangle.
 \end{split}
\end{equation}
\end{definition}
Simple calculations show that the partial pairings involved in
the above formula can be expressed as follows:
\begin{align}
\label{Eq:FiEq}
   &\langle \varphi_{2k+1},\ch_\bullet(a_s)\rangle=\frac{k!}{(2k+1)!} \int 
   \big( a_s^{-1}da_s \big)^{2k+1},\\
\label{Eq:SecEq}
   \begin{split}
   &\langle \psi_{2k},\slch_\bullet( \sigma (a_s) ,\sigma( \dot a_s ))\rangle\\
    &\quad=
   \frac{k!}{(2k)!} \int' \! (\sigma (a_s)^{-1}\sigma ( \dot a_s) ) 
   \big( (\sigma (a_s))^{-1}d (\sigma(a_s))\big)^{2k}.
   \end{split}
\end{align}
By Theorem.~\ref{Thm:NatChern}, the divisor flow is a homotopy invariant, 
which immediately entails the following result.
\begin{theorem}
\label{Thm:HomInvDF}
  Assume that $\sigma :\cA \rightarrow \cB$ is a continuous 
  surjective unital homomorphism of Fr\'echet algebras, with $\cA$ good,
  and let $C$ be an odd relative cycle 
  over $(\cA,\cB)$ with character $(\varphi_{2k+1},\psi_{2k}) :=\character C$.
  Then the divisor flow with respect to $C$ is
  well-defined on 
  $\pi_1 \big(\Ell_\infty (\cA), \GL_\infty (\cA)\big)$.
  Moreover, the divisor flow map 
  \[ \DF_C: \pi_1 \big(\Ell_\infty (\cA), \GL_\infty (\cA)\big)
  \rightarrow \mathbb{C} \]
  is additive with respect to composition of paths. 
  Finally, if in addition $\sigma :\cA \rightarrow \cB$ is 
  a morphism of local Banach algebras, then the divisor flow defines a 
  homomorphism 
  \[ \DF_C: K_1(\cA,\cB) \rightarrow \mathbb{C} . \]
\end{theorem}

\subsection{Divisor flow associated to an even relative cycle}
\label{evencase}

In this section we treat the even case.  
In many respects it is parallel to the odd case. 
Since the relative $K_0$-theory of algebras is quite well-covered in the literature,
cf. in particular ~\cite[Theorem 5.4.2]{Bla:TOA} and
{\textsc{Higson--Roe} \cite[Sec. 4.3]{HigRoe:AKH}, 
we will be rather brief and state the results without proof.

Throughout the entire section $\cA$ and $\cB$ will denote
two unital Fr\'echet algebras, and 
  \begin{equation}\label{ExSeqAlg}
    0 \longrightarrow \cJ \longrightarrow  \cA 
    \overset{\sigma}{\longrightarrow} \cB \longrightarrow 0
  \end{equation}
is an exact sequence of Fr\'echet algebras and unital homomorphisms, with
$\cA$ (and hence $\cB$) a good Fr\'echet algebra. 

We denote by $\Proj_\infty(\cB)$ the set of idempotents in the stable
matrix algebra $\mathfrak{M}_\infty(\cB)$ and by 
\begin{equation}
   \AP_\infty(\cA):=
   \bigl\{ x\in\mathfrak{M}_\infty(\cA)\,\big|\, 
   \sigma(x)\in\Proj_\infty(\cB)\bigr\}
\end{equation}
the set of \emph{almost idempotents} in $\cA$.   

The relative $K_0$-group $K_0(\cA, \cB)$ can be described as follows:
let $V(\cA,\cB)$ be the set of triples $(e,f,\gamma)$, where
$e,f\in \Proj_\infty(\cA)$ and $\gamma:[0,1]\to\Proj_\infty(\cB)$ is
a continuous path with $\gamma(0)=\sigma(e), \gamma(1)=\sigma(f)$. 
A triple $(e,f,\gamma)$
is \emph{degenerate}, if $\gamma$ has a continuous lift
$\widetilde\gamma:[0,1]\longrightarrow \Proj_\infty(\cA)$ with 
$\sigma\circ\widetilde \gamma=\gamma$ and $\widetilde\gamma(0)=e, \widetilde\gamma(1)=f$. 
A homotopy of triples
is a continuous path $(e_t,f_t,\gamma_t)$ of triples. Triples
can be added in the obvious way. $K_0(\cA,\cB)$ is the group
obtained from $(V(\cA,\cB),\oplus)$ by identifying homotopic triples
and dividing by the submonoid of homotopy classes of degenerate triples,
cf.~\cite[Def. 4.3.3]{HigRoe:AKH}.
Note that \cite{HigRoe:AKH} defines the triples differently by considering 
$(e,f,z)$, where $z$ is a Murray--von Neumann equivalence between $\sigma(e)$ 
and $\sigma(f)$. 
The resulting $K$-groups are the same since in the stable algebra
two projections are Murray--von Neumann equivalent iff they are homotopic
\cite[Sec. 4.2--4.4]{Bla:TOA}.

By strong excision in topological $K$-theory \cite[Theorem 4.3.8]{HigRoe:AKH} 
the natural map
\begin{equation}
   K_0(J)=K_0(\widetilde J,\widetilde J/J)\longrightarrow K_0(A,B),
   \quad [e]-[f]\mapsto [(f,e,1)]
\end{equation}
is an isomorphism. Note that elements of $K_0(J)$ are represented in
the form $[e]-[f]$ with $e,f\in\Proj_\infty(\widetilde J)$ such that
$\sigma(e)=\sigma(f)$. Since $\cA$ is unital, $\widetilde J$ is naturally
a subalgebra of $\cA$.

We now proceed to give an alternative description of the relative 
$K_0$-group, in the spirit of the discussion preceding
Theorem \plref{Thm:RelKEll}. 
Let $\Omega(\AP_\infty(\cA),\Proj_\infty(\cA))$ be the set of continuous paths
in $\AP_\infty(\cA)$ with endpoints in $\Proj_\infty(\cA)$.
Direct sum makes $\Omega\approj$ and also the homotopy set $\pi_1\approj$ 
into a monoid. A path $\gamma\in \Omega(\AP_\infty(\cA),\Proj_\infty(\cA))$
is called {\it degenerate}, if $\gamma$ maps into $\Proj_\infty(\cA)$. 
The quotient of $\pi_1\approj$ by the submonoid of homotopy classes of 
degenerate paths,  or in other words of homotopy classes of constant paths, 
is a group, which will be denoted $\widetilde\pi_1\approj$.
There is an obvious homomorphism 
\begin{equation}
\widetilde\pi_1(\AP_\infty(\cA),\Proj_\infty(\cA))\longrightarrow K_0(\cA,\cB),
\quad \gamma\mapsto (\gamma(0),\gamma(1),\sigma\circ\gamma) ,
\end{equation}
that is easily seen to be an isomorphism. Using excision we obtain
the analogue of Theorem \plref{Thm:RelKEll}, (2) in the even case.

\begin{theorem}\label{Thm:RelKElleven}
The canonical homomorphism  
\begin{displaymath}
  \begin{split}
    \kappa: \: K_0 (J) \, & \rightarrow    \widetilde\pi_1 
    \big( \AP_\infty (\cA) , \Proj_\infty (\cA)\big), \\
    [e]-[f] \, & \mapsto \big[ [0,1] \ni s \mapsto  (1-s)f + s e \in  
    \AP_\infty(\cA) \big]
  \end{split}
\end{displaymath} 
is an isomorphism.
\end{theorem}

Turning next to the even Chern character,  we recall that  the Chern character 
of an idempotent
$e\in\Proj_\infty(\cA)$ is given by the formula
\begin{equation}
  \ch_\bullet(e) := 1+ \sum_{k=1}^\infty (-1)^k \frac{(2k)!}{k!}  
  \tr_{2k} \Big( \big(e - \frac{1}{2} \big)\otimes e^{\otimes (2k)}\Big)
  \in HP_0 (A).
\end{equation}
If $(e_s)_{0\leq s \leq 1}$ is a smooth path of idempotents, then the 
transgression formula reads
\begin{equation}   
\label{evenh}
  \frac{d}{d s}\ch_\bullet(e_s)
  = (b+B)\slch_\bullet(e_s,(2e_s-1)\dot e_s);
\end{equation}
here the secondary Chern character $\slch_\bullet$ is given by
\begin{equation}
  \slch_\bullet(e,h):=\iota (h)\ch_\bullet(e),
\end{equation}
where the map $\iota (h)$ is defined by
\begin{equation}\begin{split}
 \iota (h)&(a_0\otimes a_1\otimes \ldots\otimes a_l)\\
        &=\sum_{i=0}^l (-1)^i 
        (a_0\otimes \ldots\otimes a_i\otimes h\otimes a_{i+1}
        \otimes \ldots\otimes a_l).
\end{split}
\end{equation}

In analogy to Proposition \plref{Prop:TransRelCyc} we now define
the Chern character of a path $(f_s)_{0\le s\le 1}\in\Omega\approj$
by putting
\begin{equation}
\begin{split}
   \ch_\bullet&\big((f_s)_{0\le s\le 1}\big)\\
   &:=\Big(\ch_\bullet(f_1)-\ch_\bullet(f_0),
   -\int_0^1 \slch_\bullet \big(\sigma (f_s) , 
   \sigma ((2f_s-1)\dot f_s )\big)   \, ds \Big).
\end{split}
\end{equation}
By the transgression formula $\ch_\bullet \big( (f_s)_{0\le s\le 1} \big)$ 
indeed is a relative cyclic cycle.

The even counterpart of the secondary transgression 
formula for $\slch_\bullet$ has already appeared in
\textsc{Moscovici--Wu}~\cite[Lemma 1.11]{MosWu:ITS},
and can be stated as follows.
Let $e_{s,t}\in\Proj_\infty(\cA)$, $s,t\in [0,1]$ be a smooth two parameter
family of idempotents. Put
\begin{equation}
    a_t:=(2e-1)\partial_t e,\quad a_s:=(2e-1)\partial_s e.
\end{equation}
Then
\begin{equation}
\label{Eq:SlChTransgressEven}
  \frac{\partial}{\partial s} \slch_\bullet (e , a_t ) 
  - \frac{\partial}{\partial t} \slch_\bullet (e , a_s) 
  = (b+B) \,  \sllch_\bullet (e, a_s, a_t),
\end{equation}
where
\begin{equation}
  \begin{split}
    \sllch_\bullet&(e,h_1,h_2)\\&:=\iota(h_1)\iota(h_2)
    \ch_\bullet(e)+R_\bullet([h_1,h_2],e)-2Q_\bullet([h_1,h_2],e).
  \end{split}
\end{equation}
We omit reproducing the explicit
formulas for $R$ and $Q$, which are quite involved,
and refer instead to ~\cite[Prop. 1.14]{MosWu:ITS}. 

In view of the above secondary transgression formula, $\ch_\bullet$ is seen to
be homotopy invariant. Since degenerate paths are null-homotopic and 
$\ch_\bullet$ obviously vanishes on constant paths, $\ch_\bullet$ does 
descend to a map defined on $\widetilde\pi_1\approj\cong K_0(\cA,\cB)$.
Thus $\ch_\bullet\circ \kappa$ coincides with the standard Chern character 
on $K_0(J)\cong K_0(\cA,\cB)$, giving the even case of 
Theorem \plref{Thm:NatChern}.

Finally, we are going to define the analogue of the divisor flow in the even 
case. 
As in Section \ref{sec1.3}, we fix a relative cycle
$C$ over $(\cA,\cB)$, this time of degree $2k$,  with
character $\character C=(\varphi_{2k},\psi_{2k-1})$.

Let $(f_s)_{0\le s\le 1}$ be a smooth path in $\AP_\infty(\cA)$ with
endpoints in $\Proj_\infty(\cA)$. Then the relative pairing between  
$\ch_\bullet \big( (f_s)_{0\leq s \leq 1} \big)$ and 
the character $(\varphi_{2k},\psi_{2k-1})$ of $C$ defines the
{\it even divisor flow}
\newcommand{\constevendiv}{\frac{(-1)^{k+1}}{(2\pi i)^k}}
\begin{equation}\label{Eq:dfeven}
 \begin{split}
    \DF_C \big(  &(f_s)_{0\leq s\leq 1} \big) := 
    \DF \big( (f_s)_{0\leq s\leq 1} \big) := \\
    := \, &  \constevendiv \,
    \big\langle \character C,
    \ch_\bullet ((f_s)_{0\le s\le 1})\big\rangle \\
    = \, &\constevendiv \,
    \Big( \big\langle 
    \varphi_{2k},\ch_\bullet(f_1)\big\rangle -
    \big\langle \varphi_{2k}, \ch_\bullet(f_0)\big\rangle\Big) \\  
    & +\frac{(-1)^{k}}{(2\pi i)^k} \,
    \big\langle \psi_{2k-1},\int_0^1 
    \slch_\bullet(\sigma (f_s), \sigma ((2f_s-1)\dot f_s))ds \big\rangle.
 \end{split}
\end{equation}
Explicitly, the partial pairings entering in the above formula are given by
\begin{align}
\label{Eq:FiEqeven}
   &\langle \varphi_{2k},\ch_\bullet(f_s)\rangle=
   \frac{(-1)^k}{k!} 
   \int \big( f_s-\frac 12 \big)(df_s)^{2k},\\
\label{Eq:SecEqeven}
   \begin{split}
   &\langle \psi_{2k-1},\slch_\bullet( \sigma (f_s) ,
   \sigma( (2f_s-1)\dot f_s ))\rangle = \\
    &\quad=
  \frac{(-1)^k}{ (k-1)!} \int' \! \sigma ( (2f_s-1)\dot f_s)
   \big(d (\sigma(f_s))\big)^{2k-1}.
   \end{split}
\end{align}
Under the same assumptions on $\cA$ and 
$\cB$ as in Theorem \plref{Thm:HomInvDF},
the even divisor flow then defines a homomorphism
\begin{equation}
  \DF_C: K_0(\cA,\cB) \rightarrow \mathbb{C} . 
\end{equation}
\section{Divisor flow on suspended 
         pseudodifferential operators}
\subsection{Pseudodifferential operators with parameter and 
regularized traces}
\label{Sec:Parametric}
We start by establishing some specific notation and recalling a few basic notions.
Let $U\subset \R^n$ be an open
subset. We denote by $\sym^m(U;\R^N)$, $m\in \R$, the space of symbols 
of H\"ormander type $(1,0)$ (\textsc{H\"ormander} \cite{Hor:FIOI}, 
\textsc{Grigis--Sj{\o}strand} \cite{GriSjo:MAD}). 
More precisely, $\sym^m(U;\R^N)$ consists of those 
$a\in \cC^\infty(U\times \R^N)$ such that for multi-indices 
$\alpha\in \Z_+^n,\gamma\in \Z_+^N$ and each compact subset $K\subset U$ 
we have an estimate
\begin{equation}\label{ML-G2.1}
    \bigl|\partial_x^\alpha\partial_\xi^\gamma a(x,\xi)\bigr|
  \le C_{\alpha,\gamma,K} (1+|\xi|)^{m-|\gamma|}, \quad x\in K, \xi\in \R^N.
\end{equation}
The best constants in \eqref{ML-G2.1} provide a set of 
semi-norms which endow
$\sym^\infty (U;\R^N):=\bigcup_{m\in\R}S^m(U;\R^N)$ with the structure of a
Fr{\'e}chet-algebra. 

A symbol $a\in S^m(U;\R^N)$ is called
\emph{$\log$-polyhomogeneous} (cf.~\textsc{Lesch} \cite{Les:NRP}), if it has an 
asymptotic expansion in $\sym^\infty (U;\R^N)$ of the form
\begin{equation}\label{ML-G2.2}
    a\sim\sum\limits_{j=0}^\infty a_j \quad   
    \text{ with $a_j=\sum_{l=0}^{k_j} b_{lj}$}, 
   \end{equation}
where $a_j\in \cC^\infty(U\times \R^N)$ and 
$b_{lj}(x,\xi)=\tilde b_{lj}(x,\xi/|\xi|)|\xi|^{m_j}\log^l|\xi|$ for
$\xi\ge 1$. Here, $(m_j)_{j\in \N}$ is a decreasing sequence with
$m\ge m_j\searrow -\infty$.

A symbol $a\in \sym^m(U;\R^N)$, $m\in \R$, is called \emph{classical}, 
if it is $\log$-poly\-homo\-geneous
and $m_j=m-j$, $k_j=0$ for all $j$. In other words, $\log$-terms do not
occur in the expansion \eqref{ML-G2.2} and the degrees of
homogeneity in the asymptotic expansion decrease in
integer steps. 

The space of $\log$-polyhomogeneous symbols with $m=m_0$ in the expansion 
\eqref{ML-G2.2} is denoted by $\PS^m(U;\R^N)$ and the space of classical
symbols in $\sym^m(U;\R^N)$ by $\CS^m(U;\R^N)$.
Note the little subtlety that in view of the $\log$-terms
we only have 
$\PS^m\subset\bigcap\limits_{\varepsilon > 0} \sym^{m+\varepsilon}$
but $\PS^m\not\subset \sym^m$.

Fix $a\in\sym^m(U;\R^n\times\R^p)$ (resp.~$\in \CS^m(U;\R^n\times\R^p)$).
For each fixed
$\mu_0$ we have $a(\cdot, \cdot, \mu_0) \in \sym^m (U; \R^n)$
(resp.~$\in\CS^m(U;\R^n))$ and hence
we obtain a family of pseudodifferential operators
parameterized over $\R^p$ by putting
\begin{equation}
\begin{split}
 \big[ \Op( &a(\mu_0) ) \, u \big] \, (x):=  \big[ A(\mu_0) \, u \big] (x)\\
      &:= \int_{\R^n} \, e^{i \langle x,\xi \rangle} \, 
      a(x,\xi,\mu_0) \, \hat{u} (\xi ) \, \dbar \xi \\
      &= \int_{\R^n}\int_U \, e^{i \langle x-y,\xi \rangle} \, 
      a(x,\xi,\mu_0) \, u(y) \dbar y \dbar \xi .
\end{split}
\end{equation}
From now on, $\dbar$ denotes $(2\pi)^{- n/2}$-times the Lebesgue measure 
on $\R^n$.  

We denote by $\pdo^m(U;\R^p)$ the set of all $\Op(a)$
(plus parameter dependent smoothing operators), where
$a\in \sym^m(U;\R^p)$, and by $\CL^m(U;\R^p)$ the corresponding
algebra of classical parameter dependent pseudodifferential
operators.

\begin{remark}
In case $p = 0$ we obtain the usual (classical)
pseudodifferential operators of order $m$ on $U$. 
Parameter dependent pseudodifferential operators play a crucial
role, e.g., in the construction of the resolvent expansion
of an elliptic operator (\textsc{Gilkey} \cite{Gil:ITHEASIT}). The definition
of the parameter dependent calculus is not uniform in the literature.
It will be crucial in the sequel that differentiating by the parameter
reduces the order of the operator. This is the convention e.g.
of \cite{Gil:ITHEASIT} but differs from the one in 
\textsc{Shubin} \cite{Shu:POST}. 
In \textsc{Lesch--Pflaum} \cite[Sec.~3]{LesPfl:TAP} 
it is shown that parameter dependent pseudodifferential operators can
be viewed as translation invariant pseudodifferential
operators on $U\times \R^p$ and therefore our convention
of the parameter dependent calculus captures Melrose's
suspended algebra from \cite{Mel:EIF}. 

For a smooth manifold $M$ and a vector bundle $E$ over $M$ we define the
space $\CL^m (M,E; \R^p)$ of classical parameter dependent pseudodifferential
operators between sections of $E$ in the usual way by patching together local
data.

The \emph{suspension} of a (Fr{\'e}chet) algebra $\cA$ is, by definition,
 the algebra
$C_0(\R,\cA)$ of all continuous functions on $\R$ with values in 
$\cA$ which vanish at infinity. However, in the pseudodifferential world
the `right' setting for considering
pseudodifferential operator valued functions is 
the parameter dependent calculus. In this sense, $\CL^0(M,E;\R)$
could be viewed as the \emph{pseudodifferential suspension} (with unit) 
of the algebra $\CL^0(M,E)$.
Indeed, it was shown by \textsc{S. Moroianu}~\cite{Mor:TSPO} that
the $K$-theory of $\CL^0(M,E;\R)$ does 
coincide with the $K$-theory of $C_0(\R,\CL^0(M,E))^{++}$ (continuous
functions $f: \R \rightarrow \CL^0(M,E)$ admitting a limit at $\pm\infty$
proportional to the identity). 

As a side remark, we should mention that 
there remain some unsettled questions regarding this picture,
which deserve further attention. The true
suspension is not contained  
in $\CL^0(M,E;\R)$, nor vice versa. 
The above mentioned $K$-theoretical isomorphism is abstract and
not obviously canonical, although the latter is very
conceivable. It would
be helpful to know that an element $A\in \CL^0(M,E;\R)$ has limits
(in the operator norm on $L^2$-space) as $\mu\to\pm\infty$. 
However, in this respect another disconcerting thing happens. Since the (ordinary
pseudodifferential) leading symbol of $A(\mu)$ is independent
of $\mu$, it is clear that $\lim\limits_{\mu\to\infty} A(\mu)
\equiv A(0)$ modulo lower order terms. On the other hand, at the
$K$-theoretical level, $\CL^0(M,E;\R)$ is the algebra of pseudodifferential
operator valued functions whose limits are proportional
to the identity.  
\end{remark}

We now fix a compact smooth manifold $M$ without boundary of dimension $n$.  
Denote the coordinates in $\R^p$ by $\mu_1,\ldots,\mu_p$ and
let $\polyn$ be the algebra of polynomials
in $\mu_1,\ldots,\mu_p$. By slight abuse of notation we denote
by $\mu_j$ also the operator of multiplication by the
$j$-the coordinate function. Then we have
\begin{equation}
\begin{split}
   &\partial_j:\CL^m(M,E;\R^p)\rightarrow \CL^{m-1}(M,E;\R^p),\\
   &\mu_j:\CL^m(M,E;\R^p)\rightarrow  \CL^{m+1}(M,E;\R^p).
\end{split}
\end{equation}
Also $\partial_j$ and $\mu_j$ act naturally on $\PS^\infty (\R^p)$
and $\polyn$ and hence on the quotient
$\PS^\infty (\R^p)/\polyn$.
After these preparations we can summarize some of the results
of \cite{LesPfl:TAP}. 

Let $E$ be a smooth vector bundle on $M$ and consider $A\in\CL^m(M,E;\R^p)$
with $m+n< 0$. Then for $\mu\in\R^p$
the operator $A(\mu)$ is trace class hence we may define the function
$\TR(A):\mu\mapsto \tr(A(\mu))$. The map $\TR$ is obviously tracial, 
i.e.~$\TR(AB)=\TR(BA)$, and commutes with $\partial_j$ and $\mu_j$. 
In fact, the following theorem holds.

\begin{theorem}
\textup{\cite[Theorems 2.2, 4.6 and Lemma 5.1]{LesPfl:TAP}}
There is a unique linear extension
\[\TR:\CL^\infty (M,E;\R^p)\rightarrow \PS^\infty (\R^p)/\polyn\]
of $\TR$ to operators of all orders such that
\begin{enumerate}
\item \label{thm21.1}$\TR(AB)=\TR(BA)$, i.e. $\TR$ is tracial.
\item \label{thm21.2}$\TR(\partial_j A)=\partial_j \TR(A)$ for $j=1,...,p$.
\end{enumerate}
This unique extension $\TR$ satisfies furthermore:
\begin{enumerate}
\setcounter{enumi}{2}
\item \label{thm21.3}$\TR(\mu_j A)=\mu_j\TR(A)$ for $j=1,...,p$.
\item \label{thm21.4}$\TR(\CL^m(M,E;\R^p))\subset \PS^{m+p}(\R^p)/\polyn$.
\end{enumerate}
\end{theorem}
\begin{proof}[Sketch of Proof]
We briefly present two arguments which help explain why this 
theorem is true.

\subsubsection*{1.~Taylor expansion} 
 Given $A \in \CL^m (M,E;\R^p)$. Since differentiation by the parameter
reduces the order of the operator, the Taylor expansion around $0$
yields for $\mu\in\R^p$ (cf.~\cite[Prop.~4.9]{LesPfl:TAP})
 \begin{equation} \label{LesPfl:G1-3.7}
   A(\mu) - \sum_{|\alpha| \leq N-1} \, 
   \frac{(\partial_\mu^\alpha A)(0)}{\alpha !} \, 
   \mu^\alpha \in \CL^{m-N} (M,E).
 \end{equation}
Hence, if $N$ is so large that $m-N+n < 0$, then the difference
\eqref{LesPfl:G1-3.7} is trace class and we put
 \begin{equation}
 \label{LesPfl:G1-3.8}\begin{split}
   \TR(A)(\mu) 
    := \tr &\Big( A(\mu) - \sum_{|\alpha| \leq N-1}
   \frac{(\partial_\mu^\alpha A) (0)}{\alpha !} \, \mu^\alpha
   \Big) \; \mod {\polyn}.
 \end{split}
\end{equation}
Since we mod out by polynomials, the result is in fact independent
of $N$. This defines $\TR$ for operators of all order and the properties
(1)--(3) are straightforward to verify.
However,  \eqref{LesPfl:G1-3.7} does not give any asymptotic 
information and hence does not justify the fact that $\TR$ takes
values in $\PS$.

\subsubsection*{2.~Differentiation and integration} 
Given $A\in \CL^m (M,E;\R^p)$, then 
$\partial^\alpha A\in \CL^{m-|\alpha|}(M,E,\R^p)$ which again is of
trace class if $m-|\alpha|+ \dim M <0$. Now integrate the function 
$\TR(\partial^\alpha A)(\mu)$
back. Since we mod out polynomials this procedure is independent
of $\alpha$ and the choice of antiderivatives. This integration procedure
also explains the possible occurrence of $\log$-terms in the
asymptotic expansion and hence why $\TR$ ultimately takes values in
$\PS$. For details, see \cite[Sec.~4]{LesPfl:TAP}. 
\end{proof}

\medskip
Composing any linear functional on $\PS^\infty (\R^p)/\polyn$
with $\TR$ yields a trace on $\CL^\infty (M,E;\R^p)$. The
\emph{regularized trace} is obtained by composition
with the following regularization of the multiple
integral. If $f\in\PS^m(\R^p)$ then
\begin{equation}
    \int_{|\mu|\le R} f(\mu)d\mu\sim_{R\to\infty}
       \sum_{\alpha\to-\infty}
       p_\alpha(\log R) R^\alpha, \label{G2-4.1}
\end{equation}
where $p_\alpha$ is a polynomial of degree $k(\alpha)$. The 
\emph{regularized integral}
$\displaystyle \regint_{\R^p} f(\mu) d\mu$ is, by definition, the 
constant term in this asymptotic expansion. This 
notion of regularized integral is close to 
the Hadamard \emph{partie finie}. It has a couple of peculiar
properties, cf.~\cite{Mel:EIF}, which were further investigated in 
\cite[Sec.~5]{Les:NRP} and \cite{LesPfl:TAP}. 
The most important features are a modified change of
variables rule for linear coordinate changes
and, as a consequence, the fact that Stokes' theorem does
not hold in general. The failure of Stokes' theorem gives rise,
as we shall see below, to a `symbolic' \emph{formal trace}.

Both the regularized and the formal trace are best explained
in terms of differential forms on $\R^p$ with coefficients
in $\CL^\infty (M,E;\R^p)$. Let 
$\Lambda^\bullet:= \Lambda^\bullet (\R^p)^*=\C[d\mu_1,\ldots,d\mu_p]$ 
be the exterior algebra of the vector space $(\R^p)^*$. We put 
\begin{equation}\label{ML-G2.8}
    \Omega_p:=\CL^\infty (M,E;\R^p)\otimes \Lambda^\bullet .
\end{equation}
Then, $\Omega_p$ just consists of pseudodifferential operator-valued
differential forms, the coefficients of $d\mu_I$ being
elements of $\CL^\infty (M,E;\R^p)$.

For a $p$-form $A(\mu)d\mu_1\wedge\ldots\wedge d\mu_p$
we define the \emph{regularized trace} by
\begin{equation}\label{ML-G2.9}
     \fTR(A(\mu)d\mu_1\wedge\ldots\wedge d\mu_p)
   := \regint_{\R^p} \TR(A)(\mu)d\mu_1\wedge\ldots\wedge d\mu_p.
\end{equation}
On forms of degree less than $p$ the regularized trace
is defined to be $0$. $\fTR$ is a graded trace on the
differential algebra $(\Omega_p,\, d)$. In general, $\fTR$ is
not closed. However, its boundary, 
$$
  \lTR:= d\fTR := \fTR\circ d \, ,
$$ 
called the \emph{formal trace},
is a closed graded trace of degree $p-1$. It is
shown in \cite[Prop.~5.8]{LesPfl:TAP}, \cite[Prop.~6]{Mel:EIF}
that $\lTR$ is \emph{symbolic}, i.e.~it descends to a well-defined
closed graded trace of degree $p-1$ on
\begin{equation}\label{ML-G2.10}
  \partial \Omega_p:= \CL^\infty (M,E;\R^p)/\CL^{-\infty}(M,E;\R^p)
  \otimes\Lambda^\bullet.
\end{equation}
Denoting by $r$ the quotient map $\Omega_p\to\partial \Omega_p$
we see that Stokes' formula with `boundary' 
\begin{equation}\label{ML-G2.11}
   \fTR(d\omega)=\lTR(r\omega)
 \end{equation}
now holds by construction for any $\omega\in\Omega$.
In sum, we have constructed a cycle with boundary 
\begin{equation}
\label{regcycle}
 C^p_\text{\rm\tiny reg} :=
 (\Omega_p,\, \partial \Omega_p,\, r,\, \varrho,\, \fTR,\, \lTR)
\end{equation} 
in the sense of Definition \ref{defcyclewithboundary}.

\begin{remark}
To shorten the notation, in what follows we suppress  the subscript $_p$ in the above
cycle, and use the abbreviation 
$\csym^\infty (M,E;\R^p):=\CL^\infty (M,E;\R^p)/\CL^{-\infty}(M,E;\R^p)$. Moreover,
if no confusion can arise, we will often write $\CL^m_p$ for  
$\CL^m (M,E;\R^p)$,  $\CS^m_p$ for  
$\CS^m (M,E;\R^p)$, and finally $\csym^m_p$ for $\csym^m (M,E;\R^p)$. 
\end{remark}

\subsection{The Fr\'echet structure of the suspended algebra}
For each fixed $m\in \Z$, the space $\CL^m (M,E;\R^p)$, $p\in \N$, carries 
a natural 
Fr\'echet  topology, which will be described below. To keep the notation simple,
we assume that the vector bundle $E$ is trivial.
For the construction of the Fr\'echet structure on  $\CL^m_p$ it will be convenient
to use the global symbol calculus for pseudodifferential operators 
developed by \textsc{Widom}~\cite{Wid:CSCPO} (see also \cite{Pfl:NSRM}). 
This requires to fix a Riemannian metric on $M$. Let $d$ be the geodesic 
distance with respect to this Riemannian metric, and $\varepsilon$  the 
corresponding injectivity radius. Then choose a smooth function
$\chi : \R M \rightarrow [0,1]$ such that
$\chi (s) = 0$ for $s \geq \varepsilon$ and $\chi (s) = 1$ for 
$s\leq \frac 34 \varepsilon$.  
Put $\psi (x, y ) = \chi (d^2(x,y))$ for $x,y \in M$ and 
$\psi (v) = \chi ( \pi(v), \exp v ) $ for 
$v \in TM$, where $\pi$ is the projection of
the (co)tangent bundle. Then $ \psi$ is a cut-off function around the diagonal
of $M\times M$ resp.~around the zero-section of $TM$.  
With these data let us define two maps,
namely the {\it symbol map} $ \symb: \CL^m_p \rightarrow \CS^m_p$
and the {\it operator map} $\Op: \CS^m_p \rightarrow \CL^m_p $,  by 
putting for $A\in \CL^m_p$, $a \in \CS^m_p$ and $\mu \in \R^p$:
\begin{displaymath}
\begin{split}
  \symb(A) :&\, T^*M \times \R^p \rightarrow \C, \:
  (\xi,\mu) \mapsto \big[ A(\mu)  \psi (  . , \pi (\xi)) \, 
  e^{i \langle \exp^{-1}_{\pi (\xi)} (.), \xi  \rangle} \big] (\pi(\xi) ) , \\
  \Op (a ) (\mu) u  :&\, M \rightarrow  \C, \: 
  x \mapsto \int_{T^*_x M} \int_{T_xM} 
  a (\xi, \mu) \,  
  e^{-i\langle v , \xi \rangle} \, 
  \psi (v) \, u(\exp v) \, \dbar v \, \dbar \xi  .
\end{split}
\end{displaymath}
The maps $\Op$ and $\symb $ are quasi-inverse
to each other, which means that 
\begin{equation}\label{Eq:symbolsection}
 A -\Op (\symb (A)) \in \CL^{-\infty}_p \: \text{ and  } \:
   a -\symb (\Op (a)) \in \CS^{-\infty}_p
 \end{equation}
for all $A\in \CL^m_p$ and $a\in \CS^m_p$. As a consequence $\symb$ induces
a natural topological isomorphism
\begin{equation}\label{Eq:symbiso}
\csym^\infty_p=\CL_p^\infty/\CL_p^{-\infty}\longrightarrow \CS_p^\infty/\CS_p^{-\infty},
\end{equation}
whose inverse is induced by $\Op$.

Now choose a finite open covering $\mathcal U = ( U_j )_{j\in J}$ of $M$
together with a subordinate partition of unity $(\varphi_j)_{j\in J}$
such that for every $j\in J$ there exists a coordinate system 
$x_{(j)} : U_j \rightarrow \R^n$. For $N\in \N$ and $A\in \CL^m_p$ then put:
\begin{displaymath}
\begin{split}
  q_{m,N} (A) := 
  \sum_{j\in J} \, \sum_{|\alpha|,|\beta|,|\gamma| \leq N}
  & \, \sup_{(x,\xi)\in T^*M \atop \mu \in \R^p} 
  (1+ |\xi|^2 + |\mu|^2)^{\frac{m-|\beta|-|\gamma|}{2}} \cdot \\
  & \cdot \varphi_j (x) \, \partial^\alpha_{x_{(j)}} \partial^\beta_{\xi_{(j)}}
  \partial^\gamma_\mu  \symb (A) (x,\xi,\mu), 
  \\
  \widetilde{q}_{m,N} (A) := 
  \sum_{j\in J} \, \sum_{|\alpha|,|\beta|,|\gamma| \leq N
  \atop 0\leq l \leq N}
  & \, \sup_{(x,\xi)\in T^*M \atop \mu \in \R^p}  (1+ |\mu|^2)^{\frac{l}{2}} \cdot \\
  & \cdot \varphi_j (x) \, \partial^\alpha_{x_{(j)}} \partial^\beta_{\xi_{(j)}}
  \partial^\gamma_\mu  K_{A(\mu)-\Op(\symb (A))(\mu)} (x,\xi) ,
\end{split}
\end{displaymath} 
where $K_B$ denotes the Schwartz kernel of a pseudodifferenetial operator $B$.
One checks immediately that the $q_{m,N}$ and $\widetilde{q}_{m,N}$ are seminorms
which turn $\CL^m_p$ into a Fr\'echet space and even a Fr\'echet algebra,
in case $m=0$. With respect to this Fr\'echet topology on $\CL^m_p$, the subspace
$\CL^{-\infty}_p$ is a closed subspace, hence the quotient space $\csym^m_p$ inherits a 
Fr\'echet topology from $\CL^m_p$, and we get an exact sequence of Fr\'echet spaces
\begin{equation}
\label{Eq:ExSeqSusAlg}
 0 \longrightarrow \CL^{-\infty}_p 
   \longrightarrow \CL^m_p \overset{\sigma}{\longrightarrow} 
   \csym^m_p \longrightarrow  0 \, .
\end{equation}
Moreover, $\CL^{-\infty}_p$ is a local Banach algebra whose closure is the 
$p$-fold suspension 
$\cC_0 (\R^p,\mathcal K)$ of the algebra of compact operators $\mathcal K$.
If $m=0$, the algebra $\CL^m_p$ is a local Banach algebra as well 
(cf.~\cite[Sec.~5]{Mor:TSPO}).  

Let us note that, in a similar way, one obtains for each $m \in \Z$
the short (and topologically split) exact sequence of Fr\'echet spaces
\begin{equation}
 0 \longrightarrow \CL^{m-1}_p 
   \longrightarrow \CL^m_p \overset{\sigma_m}{\longrightarrow} 
   \cC^\infty ( S^*_p  M ) \longrightarrow  0 \, ,
\end{equation}
where $S^*_p M := \{ (\xi,\mu) \in T^*M\times\R^p\mid |x|^2 + |\mu|^2 =1 \}$
denotes the {\it $p$-suspended cosphere bundle}, 
and $\sigma_m$ is the principal symbol map for parametric pseudodifferential 
operators of order $m$. 

\subsection{Divisor flows as a relative cyclic pairing}
\label{Sec:DivFlowRelPair}
Let us consider now the suspended algebra  
$\CL^\infty_{2k+1}$ of pseudodifferential operators on a compact 
manifold $M$ with values in a bundle $E\rightarrow M$. 
Since it gives rise to the short exact sequence of local Banach algebras 
\eqref{Eq:ExSeqSusAlg},
we are in a relative situation and can apply the abstract results of the first part
to the cycle with boundary $(\Omega,\partial \Omega,r,\varrho,\fTR,\lTR)$
defined in Section \ref{Sec:Parametric}. 

Thus, let us assume to be given a smooth family 
$A_s \in \CL^\infty_{2k+1}$, $s \in [0,1],$ of elliptic 
operators of some fixed order $m \in \N$, such that 
$A_0$ and $A_1$ are invertible. 
According to Prop.~\ref{Prop:TransRelCyc}, the family $A_s$ gives rise to the
relative cyclic cycle
$$\Big( \ch_\bullet (A_1) - \ch_\bullet (A_0) , \, - \int_0^1 
 \slch \big(\sigma (A_s) \, \sigma (\dot A_s )\big) \, ds \Big).
 $$
As explained in Section \ref{Sec:AbsDivisor}, 
one can pair this cyclic cycle with the character 
$(\varphi_{2k+1},\psi_{2k}) := \character C^{2k+1}_\text{\rm\tiny reg}$ 
of the relative cycle $(\Omega_{2k+1},\, \partial \Omega_{2k+1},\, r,\, \varrho,\, \fTR,\, \lTR)$ 
to obtain the divisor flow of the suspended algebra of 
pseudodifferential operators.  
By Eqs.~\eqref{Eq:FiEq} and \eqref{Eq:SecEq}
\begin{equation}\begin{split}
     &\langle \varphi_{2k+1},\ch_\bullet(A)\rangle=\frac{k!}{(2k+1)!} 
     \fTR \, \bigl((A^{-1}dA)^{2k+1}\bigr),\\
     &\langle \psi_{2k},\slch_\bullet( \sigma(A),\sigma(\dot A))\rangle\\
       &\qquad=
   \frac{k!}{(2k)!} \lTR \, \Big(\sigma(A)^{-1} \sigma(\dot A) 
   \big( \sigma(A)^{-1}d\sigma(A)\big)^{2k}\Big) ,
\end{split}
\end{equation}
and therefore the divisor flow has the form 
\begin{equation} \label{df}
\begin{split} 
      \DF  & \, \big( (A_s)_{0\le s\le 1} \big) = \\
      = \, & \frac{k!}{(-2\pi i)^{k+1}(2k+1)!}
      \Big(\fTR \, \big((A_1^{-1} dA_1)^{2k+1} \big)- 
       \fTR \, \big((A_0^{-1} dA_0)^{2k+1}\big) \Big) \\
      & - \frac{k!}{(-2\pi i)^{k+1}(2k)!} 
    \int_0^1 \lTR \, \Big(\sigma (A_s)^{-1} \, 
    \sigma (\partial_s A_s) \, \big(\sigma(A_s)^{-1} \, 
    d \sigma(A_s)\big)^{2k} \Big) \, ds . 
\end{split}
\end{equation}
For $k=0$ this is precisely the divisor flow originally defined by 
\textsc{Melrose}~\cite{Mel:EIF},
while for $k>0$ it gives its generalization by 
\textsc{Lesch-Pflaum}~\cite[Prop.~6.3]{LesPfl:TAP}.

It will be convenient for what follows to introduce one additional piece of notation:
for every $m \in \N$,  $\Ell^m_\infty (\CL^\infty_p)$ will denote
the space of elliptic elements of order $m$ 
in $\CL^\infty_p$, and  
$\pi_1 \big(\Ell_\infty^m (\CL^\infty_p), \GL_\infty 
(\CL^\infty_p)\big)$ stands for the fundamental groupoid of the space of elliptic 
elements  of order $m$ relative the invertible ones.  

Theorem~\ref{Thm:HomInvDF}
specializes to the present situation and reads as follows.

\begin{theorem}
  For each $m \in \N$, the odd divisor flow defines a map 
  \begin{displaymath} 
   \DF : \pi_1 \big(\Ell_\infty^m (\CL^\infty_{2k+1}), \GL_\infty 
  (\CL^\infty_{2k+1})\big) \rightarrow \C, 
  \end{displaymath}
  which is additive with respect to composition of paths. 
  Furthermore, it induces a homomorphism from 
  $K_1 (\CL^0_{2k+1} ,\csym^0_{2k+1})$ to $\C$.
  \end{theorem}

As will be seen shortly, the latter homomorphism actually establishes
the isomorphism $ K_1 (\CL^0_{2k+1} ,\csym^0_{2k+1}) \cong \Z$.
(Cf. Prop. \ref{dfintegral} below).

\begin{remark}
 It is fairly obvious  that for a fixed 
  $A \in \GL_\infty (\CL^\infty_p)$ of order $m$ the space 
  $\pi_1 \big(\Ell_\infty^m (\CL^\infty_{2k+1}), \GL_\infty 
  (\CL^\infty_{k+1});A\big)$ 
  of homotopy classes starting at $A$ is naturally isomorphic to the set
  $\pi_1 \big(\Ell_\infty (\CL^0_{2k+1}), \GL_\infty (\CL^0_{k+1});I\big)$  
  via the map 
  \[ (A_s)_{0\leq s \leq 1} \mapsto (A^{-1} A_s)_{0\leq s \leq 1}.\]
\end{remark}

\subsection{Log-additivity and integrality of the divisor flow}
Let $A, B\in\CL^p_{2k+1}$ be invertible. The expression
\begin{equation}\label{pareta}
     \fTR \big( (A^{-1}dA)^{2k+1} \big)
\end{equation}
occurring in the definition of the divisor flow has been investigated
in \cite{Mel:EIF} in the case $k=0$ and in \cite{LesPfl:TAP} in general.
In the case $k=0$,  \eqref{pareta} extends to a homomorphism
from $K_1^{\textup{alg}}(\CL^\infty_1)\to \Z$. In the case $k=1$ it
was shown in \cite[Remark 6.8]{LesPfl:TAP} that  
\begin{equation}
\label{ML-additivity2}
\begin{split}
    \TR((&(AB)^{-1}d(AB))^{3})  -\TR((A^{-1}dA)^{3}) - \TR((B^{-1}dB)^{3})\\
         &=-3d\TR(A^{-1}dA\wedge dB B^{-1}),
\end{split}
\end{equation}
and so the difference on the left hand side is symbolic.
For $k\ge 2$, it can also be shown that the difference on the left hand side
of \eqref{ML-additivity2} equals $d\TR$ of a (noncommutative) polynomial
in $A,B,dA,dB$ and hence is symbolic, too. 
This being said, the following result may come as a surprise.

\begin{theorem} 
\label{additivity} Let $A_s\in \CL^m_{2k+1}$
and $B_s\in\CL^n_{2k+1}$ with $s\in [0,1]$ and $m,n\in \Z$
be admissible paths of elliptic elements. Then we have
\begin{equation}\label{eq:additivity}
   \DF \big( (A_s B_s)_{0\le s\le 1} \big)=
   \DF \big( (A_s)_{0\le s\le 1} \big)+ 
   \DF \big( (B_s)_{0\le s\le 1} \big).
\end{equation}
\end{theorem}

 For $k=0$ the theorem is trivial, while for $k=1$ it follows
from \eqref{ML-additivity2}. The general case is more subtle though,
and its proof will
be accomplished in a series of steps.
\subsubsection*{1.~Reduction to constant path}
\begin{lemma}\label{ML-L2.7} 
\begin{enumerate}
\item If $B\in\GL_\infty \big( \CL^n_{2k+1}\big) $,
then the \emph{constant} path $B$ has vanishing divisor flow.
\item Let $(A_s)\in\CL^m_{2k+1}$ be an admissible path
of elliptic elements and let $B\in\CL^0_{2k+1}$ be
invertible. Then
\[ \DF \big( (A_sB)_{0\le s\le 1} \big)=
   \DF \big( (BA_s)_{0\le s\le 1} \big) =\DF \big( (A_s)_{0\le s\le 1} \big). 
\]
\end{enumerate}
\end{lemma}
\begin{proof} The first claim is obvious from the definition. The second
claim follows from the well-known homotopy
\begin{equation}\label{homotopy} 
      \begin{pmatrix} A_s&0\\ 0 & A_0\end{pmatrix}
      \begin{pmatrix} \cos t & \sin t\\ -\sin t &\cos t \end{pmatrix}
      \begin{pmatrix} B&0\\ 0&I\end{pmatrix}
      \begin{pmatrix} \cos t & -\sin t\\ \sin t &\cos t \end{pmatrix}
\end{equation}
which shows that $\begin{pmatrix}A_sB &0\\ 0&A_0\end{pmatrix}$
and $\begin{pmatrix} A_s &0\\ 0&A_0B
\end{pmatrix}$
are homotopic. From the homotopy invariance of the divisor flow and (1)
we infer that 
$\DF\big( (A_sB)_{0\le s\le 1} \big) =\DF\big( (A)_{0\le s\le 1}\big)$. 
The proof of $\DF\big( (BA_s)_{0\le s\le 1} \big) =
\DF\big( (A_s)_{0\le s\le 1} \big)$ is analogous.
\end{proof} 

We emphasize that this argument fails if the degree of $B$ is different from
$0$, because then \eqref{homotopy} is not a valid homotopy of admissible paths!

\begin{lemma}\label{ML-L2.9} To prove Theorem \plref{additivity}
it suffices to show that for each $n\in\R$ there exists
an invertible $B_n\in\CL^n_{2k+1}$ such that
for admissible paths $A_s\in\CL^0_{2k+1}$, $s\in [0,1]$
with $A_0=I$ one has
\begin{equation}\label{ML-G2.32}
   \DF\big( (A_s B_n)_{0\leq s \leq 1} \big)=
   \DF\big( (B_n A_s)_{0\leq s \leq 1} \big)=
   \DF\big( (A_s))_{0\leq s \leq 1} \big).
 \end{equation}
\end{lemma}
\begin{proof} In the square $[0,1] \times [0,1] $ the path $s\mapsto (s,s)$ is homotopic
to the concatenation of the paths $s\mapsto (s,0)$ and $s\mapsto (1,s)$.
Consequently the path $s\mapsto A_sB_s$ is homotopic to the concatenation
of the paths $s\mapsto A_s B_0$ and $s\mapsto A_1 B_s$. Hence we are reduced
to the case that one of the two paths is constant. 

Now suppose that the condition of the lemma is fulfilled, that 
$(A_s)_{0\le s\le 1}$ is
an admissible path of elliptic elements in $\CL^m_{2k+1}$ and that
$B'\in\CL^n_{2k+1}$ is invertible. Then applying Lemma \ref{ML-L2.7}
we find
\begin{equation}\label{ML-G2.33}
\begin{split}
      \DF&\big( (A_sB')_{0\le s\le 1} \big)=
      \DF\big( (A_sA_0^{-1} B_{m+n})_{0\le s\le 1}(B_{m+n}^{-1}A_0 B')\big)\\
      &=\DF \big( (A_sA_0^{-1} B_{m+n})_{0\le s\le 1}\big)=\DF\big( (A_sA_0^{-1})_{0\le s\le 1} \big).
\end{split}
\end{equation}
The right hand side is independent of $B'$.
Eq.~\eqref{ML-G2.33} applies in particular to $B'=I$ and hence 
$\DF\big( (A_sA_0^{-1})_{0\le s\le 1} \big)=
\DF\big( (A_s)_{0\le s\le 1} \big)$.
The proof for $\DF\big( (B'A_s)_{0\le s\le 1} \big)$ 
works exactly the same way.
\end{proof}

\subsubsection*{2.~Reduction to the finite-dimensional case} 
By Lemma \ref{ML-L2.9} and the homotopy invariance of the divisor flow we need 
to prove \eqref{ML-G2.32} for one representative $(A_s)_{0\le s\le 1}$ 
of each class in the relative $K_1$-group  
$K_1 (\CL^0_{2k+1} ,\csym^0_{2k+1}) $ and one constant invertible $B_n$ 
for each $n$. Next we are going to choose convenient $B_n$. Choose
a Riemannian metric on $M$, a hermitian metric on $E$ as well as
a metric connection on $E$. Denote by $\Delta^E$ the connection Laplacian
acting on sections of $E$. Then $\Delta^E$ is a non-negative self-adjoint
elliptic differential operator. Put
\begin{equation}
    B_n(\mu):=\left(\Delta^E+I+|\mu|^2\right)^{q/2}.
\end{equation}
Then $B_n$ is an invertible element in $\CL^n_{2k+1}$. It has the
nice property that the operators $B_n(\mu)$ commute and have a joint spectral
decomposition. Since there is an inclusion 
$I+ \cC_\text{\tiny c}^\infty(\R^{2k+1},\mathfrak{M}_\infty(\C)) 
\subset I+\CL^{-\infty}_{2k+1}$ 
which induces an isomorphism in
$K$-theory, we may represent a class in $K_1 (\CL^0_{2k+1} ,\csym^0_{2k+1})$ 
in the following form: 
\begin{equation}
     A_s(\mu)=I+sg(\mu)P,
\end{equation}
where $P$ is a spectral projection of $\Delta^E$ of sufficiently high rank and
$g:\R^{2k+1} \longrightarrow \End(\im P)$ is a smooth function with compact
support such that $I+g(\mu)$ is invertible. 

The spectral projection $P$ reduces the operator $A_s$ as well as $B_n$.
Although $P$ is not in the parametric calculus, a direct calculation of 
$\TR$ shows that the divisor flow of $(A_s)_{0\leq s\leq 1}$ equals the 
divisor flow of the \emph{finite-rank} family 
$(PA_sP)_{0\leq s\leq 1}$. Also, the divisor flow of $(A_sB)_{0\leq s\leq 1}$
equals the divisor flow of the finite-rank family 
$(PA_sBP)_{0\leq s\leq 1}$.

Hence we are reduced to prove Theorem \ref{additivity} in the
case $M=\{\textup{pt}\}$:

\subsubsection*{3.~The finite-dimensional case} 
Consider $M=\{ \text{pt} \}$. Then
\begin{equation}
  \cA^m:=\CL^m( \{ \text{pt} \} ,\C^N;\R^{2k+1})=
  \CS^m(\{\text{pt}\}\times\R^{2k+1}) \otimes\End{\C^N}
\end{equation}
is precisely the algebra of $\End \C^N$-valued symbols of 
H\"ormander type $(1,0)$.

The divisor flow makes perfectly sense for admissible paths of elliptic
elements of $\cA^\infty$. The  Theorem \ref{additivity} will now follow from 
Lemma \ref{ML-L2.9}, the following lemma, and the subsequent remark.
\begin{lemma} \label{ML-L2.10} 
Let $A_s\in\mathcal{A}^m, B_s\in\mathcal{A}^n$, $s\in [0,1]$ 
be admissible paths of elliptic elements. 
Furthermore, assume that $B_s$ is of the form
$B_s(\mu)=f_s(\mu)\otimes I_{\C^N}$, i.e.~$(B_s)_{0\leq s \leq 1}$ 
is a path of central elements of $\mathcal{A}$. 
Then Eq.~\eqref{eq:additivity} holds for $A$ and $B$.
\end{lemma}
\begin{remark} As a constant path we may e.g. choose
\begin{equation}
  B(\mu):=(1+|\mu|^2)^{q/2}\otimes I_{\C^N}.
\end{equation}
\end{remark}

\begin{proof} 
As noted at the beginning of this section, 
the case $k=0$ is trivial, so we assume $k\ge 1$ 
and abbreviate $\omega_s:=A_s^{-1}dA_s$ and $\eta_s:=B_s^{-1}dB_s$. Then
$\omega_s$ satisfies
\begin{equation}
      d\omega_s^{2l-1}=-\omega_s^{2l}.
\end{equation}
Since $B_s$ is a scalar function, we find
\begin{align}
    d\eta_s&=0,&(\omega_s+\eta_s)^{2l}&=\omega_s^{2l},\\
    \eta_s^2&=0,&(\omega_s+\eta_s)^{2l+1}&=\omega_s^{2l+1}-d(\omega_s^{2l-1}
    \wedge\eta_s),
\end{align}
and furthermore
\begin{equation}
     (A_sB_s)^{-1}d(A_sB_s)=\omega_s+\eta_s.
\end{equation}
Next we deal with the ingredients of the divisor flow.
\begin{equation}\label{ML-G2.35}
\begin{split}
 \fTR & (\omega_1^{2k+1})-\fTR(\omega_0^{2k+1})
      -\fTR((\omega_1+\eta_1)^{2k+1})+\fTR((\omega_0+\eta_0)^{2k+1})\\
      & = \fTR(\omega_1^{2k-1}\wedge\eta_1)-\lTR(\omega_0^{2k-1}\wedge\eta_0)\\
      & =\int_0^1 \frac{d}{ds}\lTR\Bigl((\sigma(A_s)^{-1}d\sigma(A_s))^{2k-1}\wedge 
      \sigma(B_s)^{-1}d\sigma(B_s)\Bigr)\,ds.
\end{split}
\end{equation}
We abbreviate $\widetilde \omega_s:=\sigma(A_s)^{-1}d\sigma(A_s),
\widetilde \eta_s:=\sigma(B_s)^{-1}d\sigma(B_s)$. Taking into account that
$B$ is central, one sees as in the proof of
\cite[Prop. 6.3]{LesPfl:TAP} that
\begin{equation}
   \begin{split}
      \frac{d}{ds}\lTR(&\widetilde \omega_s^{2k-1}\wedge\widetilde\eta_s)\\
     & =\lTR\Bigl(d\Bigl((\sigma(A_s)^{-1} \sigma(\partial_s A_s))\widetilde\omega_s^{2k-2}\wedge \eta_s\Bigr)\Bigr)\\
     &\qquad +
        \lTR\Bigl(\widetilde\omega_s^{2k-1}\wedge d\bigl(\sigma(B_s)^{-1}\sigma(\partial_s B_s)\bigr)\Bigr).
    \end{split}\label{ML-G2.33-r5.7}
  \end{equation}
The first summand on the right vanishes because $\lTR$ is a closed trace. 
For the second summand we find
\begin{equation}\label{ML-G2.34}
    \begin{split}
      \lTR\Bigl(&\widetilde\omega_s^{2k-1}\wedge d\bigl(\sigma(B_s)^{-1}\sigma(\partial_s B_s)\bigr)\Bigr)\\
      &= \lTR\Bigl(\sigma(B_s)^{-1}\sigma(\partial_s B_s)d\widetilde\omega_s^{2k-1}\Bigr)\\
      &=-\lTR\Bigl(\sigma(B_s)^{-1}\sigma(\partial_s B_s)\widetilde\omega_s^{2k}\Bigr)=0,
    \end{split}
\end{equation}
since the one-form $\omega_s$ commutes with 
$\sigma(B_s)^{-1} \sigma(\partial_s B_s)$ and the exponent $2k$ is even.

In sum, the left hand side of Eq.~\eqref{ML-G2.35} vanishes. Moreover, 
since $\eta_s^2=0$ we have $\fTR(\eta_s^{2k+1})=0$ for $k\ge 1$. 
Thus, in the expansion of the difference  
\begin{equation}
 \DF \big( (A_s B_s)_{0\le s\le 1} \big)- \DF \big( (A_s)_{0\le s\le 1} \big)
 - \DF \big( (B_s)_{0\le s\le 1} \big)
\label{ML-G2.35-r5.7}
\end{equation}
the terms involving $\fTR$ add up to $0$. 
Furthermore, since $\eta_s^2=0$ we have
\[ \lTR \, \Big(\sigma (B_s)^{-1} \, 
    \sigma (\partial_s B_s) \, \big(\sigma(B_s)^{-1} \, 
    d \sigma(B_s)\big)^{2k} \Big)=0.\]
Finally, since $(\omega_s+\eta_s)^{2k}=\omega_s^{2k}$, we have
\[\begin{split}
     \lTR \, \Big(&\sigma (A_sB_s)^{-1} \, 
    \sigma (\partial_s( A_sB_s)) \, \big(\sigma(A_sB_s)^{-1} \, 
    d \,\sigma(A_sB_s)\big)^{2k} \Big)\\
    &= \lTR \, \Big(\sigma (A_s)^{-1} \, \sigma (\partial_s A_s) \, \widetilde \omega_s^{2k} \Big)\\
  \end{split}
\]
in view of Eq.~\eqref{ML-G2.34}. This proves that in the expression for Eq.~\eqref{ML-G2.35-r5.7}
the terms involving $\lTR$ add up to $0$. The Lemma is proved.
\end{proof}

As a consequence of the additivity and of the
$K$-theoretic interpretation of the divisor
flow we can now prove its integrality. This generalizes
\cite[Prop. 8]{Mel:EIF}.

\begin{proposition}
\label{dfintegral}
The divisor flow defined on the homotopy groupoid 
$\pi_1 \big(\Ell_\infty^m (\CL^0_{2k+1}), \GL_\infty (\CL^0 _{2k+1})\big)$
assumes integer values. Moreover, it induces an isomorphism
\[
K_1(\CL^0_{2k+1},\csym^0_{2k+1}) \overset{\simeq}\longrightarrow  \Z.
\]
\end{proposition}

\begin{proof} Let $(A_s)_{0\le s\le 1}$ be an admissible path
of elliptic elements of $\CL^m_{2k+1}$. By the additivity
and by the fact that the divisor flow of a constant path vanishes
we find
\begin{equation}
    \DF \big( (A_s)_{0\le s\le 1} \big)=
    \DF \big( (A_sA_0^{-1})_{0\le s\le 1} \big)
\end{equation}
hence it suffices to prove integrality for admissible paths of 
$0^\text{th}$-order elements starting at $I$. By Theorem \ref{Thm:RelKEll}
it remains to prove integrality of the divisor flow for standard paths
in $I+\CL^{-\infty}_{2k+1}$.   

Consider an (ordinary) pseudodifferential projection $P\in \CL^0(M,E)$
of finite rank and a smooth map
\begin{equation}
  g:\R^{2k+1}\longrightarrow \GL(\im P) \text{ with } 
  \lim_{|\mu|\to\infty} g(\mu)=I.
\end{equation}
In other words, this means that $g$ is an element of the algebra obtained 
by adjoining a unit to $\cC^\infty_0(\R^{2k+1},\GL(\im P))$. 
The map
\begin{equation}
  T(\mu):=g(\mu)P+I-P
\end{equation}
is in $I+\CL^{-\infty}_{2k+1}$. Put
\begin{equation}\label{ML-G2.19}
  A_s:=(1-s)I+s T,\quad 0\le s\le 1.
\end{equation}
Then $(A_s)_{0\leq s \leq 1}$ is a smooth family of elliptic elements in 
$\CL^0_{2k+1}$ with $A_0=I$ and $A_1$ invertible.

Every element of $K_1(\CL^0_{2k+1},\csym^0_{2k+1})$ can be
represented in this way. Indeed, by excision and the well-known fact
that dense subalgebras which are stable under holomorphic functional
calculus have the same $K$-theory as the original algebra
\begin{equation}
  \begin{split}
     K_1(\CL^0_{2k+1},\csym^0_{2k+1}) & =K_1(\CL^{-\infty}_{2k+1})
     =K_1(\cC_0^\infty(\R^{2k+1}, \mathfrak{M}_\infty(\C))) \\
     & =K_1(\cC_0^\infty(\R^{2k+1},\C)) \cong \Z.
  \end{split}
\end{equation}
To achieve the proof it will suffice to show that the divisor flow of $(A_s)$ is integral
and that there exists a path $(A_s)$ of the form \eqref{ML-G2.19} and 
of divisor flow $1$.

Since $\sigma(A_s)=1$ and $A_0=I$, the divisor flow equals
\begin{equation}\begin{split}
   \DF\big( &(A_s)_{0\le s\le 1} \big)=(-2\pi i)^{-(k+1)}\langle 
   \varphi_{2k+1},\ch_\bullet (A_1)\rangle
   \\
   &=\frac{k!}{(-2\pi i)^{k+1}(2k+1)!}\int_{\R^{2k+1}} 
   \tr(g(\mu)^{-1}dg(\mu))^{2k+1}d\mu.
 \end{split}
\end{equation} 
 The latter is precisely the odd Chern character of $g$ 
 \cite[Prop. 1.4]{Get:OCC}.
 The odd Chern character is known to be an isomorphism from 
 $K_1(\cC_0(\R^{2k+1}))$ onto $\Z$, and so we reach the desired conclusion.
\end{proof}
\subsection{Compatibility with Bott periodicity}
Recall the exact sequence \eqref{Eq:ExSeqSusAlg} of Fr\'echet spaces
\begin{equation*}
 0 \longrightarrow \CL^{-\infty}_p 
   \longrightarrow \CL^m_p \overset{\sigma}{\longrightarrow} 
   \csym^m_p \longrightarrow  0 \, .
\end{equation*}
As has been mentioned, $\CL_p^{-\infty}$ is nothing
but the $p$-fold smooth suspension of the algebra $\CL^{-\infty}(M,E)$ of 
smoothing operators acting on sections of $E$. In turn, $\CL^{-\infty}(M,E)$ is 
a local Banach algebra whose closure is the algebra of compact operators, hence
the $K$-groups of $\CL^{-\infty}(M,E)$ are naturally isomorphic to those
of $\C$. For $m=0$ the above exact sequence consists of algebras, hence
by excision,
\begin{equation} \label{bott}
      K_i ( \CL^0_{2k+i} ,  \, \csym^0_{2k+i}) \cong
         K_i (\CL^{-\infty}_{2k+i})  \overset{\simeq}\longrightarrow  \Z \, , \qquad i = 0, 1.
\end{equation}
The latter isomorphism is of course the classical Bott isomorphism, but it
can also be regarded as being  induced by the restriction of the divisor flow pairing.
More precisely, the divisor flow pairing acquires the following topological interpretation.

\begin{theorem} \label{bottcomp}
The divisor flow pairing with the character  
$\character C^p_\text{\rm\tiny reg}$ of the relative cycle
$(\Omega_p, \, \partial \Omega_p,\, r,\, \varrho,\, \fTR,\, \lTR)$, 
$ p = 2k+i >0$, implements the Bott isomorphism at the relative 
$K$-theory level,
\begin{displaymath}
   K_i ( \CL^0_{2k+i} , \,\csym^0_{2k+i})  
   \overset{\simeq}\longrightarrow  \Z \, , \qquad i = 0, 1 \, ,
\end{displaymath}
in a manner compatible with the Bott suspension.
\end{theorem}

\begin{proof}
The fact that \eqref{bott} is an isomorphism 
follows in the odd case from Proposition \plref{dfintegral} and its proof. 
What remains to be proved is the
compatibility with the Bott suspension, which automatically implies the even case.
To this end, we will relate the construction (cf.~\eqref{ML-G2.8}--\eqref{ML-G2.11}) 
of our relative cycle
$(\Omega,\partial \Omega,r,\varrho,\fTR,\lTR)$ 
to the work of 
\textsc{Elliott--Natsume--Nest} \cite{Ell:CCO} on the cyclic cohomology
of one-parameter crossed products. 

Assume to be given a Fr\'echet algebra 
${\cA}$ and consider the trivial $\R$-action on ${\cA}$. 
Then the smooth crossed product ${\cA}\rtimes \R$ is, via 
Fourier transform, isomorphic to the smooth suspension
${\cA}\otimes \cS(\R)=\cS(\R,{\cA})$. Here $\cS(\R)$ denotes 
the space of Schwartz functions. For the definition of the smooth suspension
see \cite[2.4]{Ell:CCO}. Note that $\cS(\R)$ is nuclear which makes
dealing with topological tensor products more convenient. The smooth
suspension is a local Banach algebra with completion 
$C_0(\R)\otimes \overline{{\cA}}$, hence we have natural isomorphisms
\begin{equation}
    K_i(\cS(\R,{\cA}))\cong 
    K_i(C_0(\R)\otimes \overline{{\cA}})\cong
    K_{i+1}(\overline{\cA})\cong K_{i+1}({\cA}).
\end{equation}
In \cite{Ell:CCO} the authors produce a natural map
\begin{equation}
     \#:HC^\bullet(\cA)\longrightarrow HC^{\bullet+1}(\cS(\R,\cA))
\end{equation}
which commutes with the periodicity operator $S$ and defines
isomorphisms
\begin{equation}
       HP^\bullet(\cA)\longrightarrow HP^{\bullet+1}(\cS(\R,\cA)).
\end{equation}
Furthermore, if 
\begin{equation}
     \beta:K_i(\cA)\longrightarrow K_{i+1}(\cS(\R,\cA))
\end{equation}
denotes the Bott suspension isomorphism, then one has for $\varphi\in Z^k_\lambda(\cA)$,
$[e]\in K_k(\cA)$
\begin{equation}
    \langle \varphi, [e]\rangle= \langle \# \varphi, \beta [e]\rangle.
\end{equation}
This means that via the natural pairing between periodic cyclic cohomology and
$K$-theory $\#$ corresponds to the Bott isomorphism.
\end{proof}

We will now describe the map $\#$ in more detail on the level of cycles
such that the relation between $\#$ and our suspended relative cycle
becomes apparent. 
Obviously, smooth suspensions correspond to crossed products with trivial 
$\R$-actions.
In our case it will be more convenient to deal with the algebra
$\cS(\R)$ with the product of pointwise multiplication. In the more
general situation of \cite{Ell:CCO} one has to deal with $\cS_*(\R)$,
that is $\cS(\R)$ equipped with the convolution product. 
This has to be taken into account when one compares the formulas in
\cite{Ell:CCO} with ours. 

We consider a locally convex cycle $(\Omega,d_\cA,\tau)$ of degree
$n$ over the Fr\'echet algebra $\cA$. That is $(\Omega,d_\cA)$
is a locally convex differential graded algebra with a continuous
closed graded trace $\tau$ of degree $n$ together with a continuous
homomorphism $\cA\longrightarrow \Omega^0$. We single out an
important example:

\begin{example}[The fundamental cycle of $\cS(\R^p)$]
Consider the algebra $\cS(\R^p)$ and put
\begin{equation}
   \cE^\bullet:= \cS(\R^p)\otimes \Lambda^\bullet,
\end{equation}
(cf.~\eqref{ML-G2.8}). This means that $\cE^\bullet$ just consists
of differential forms on $\R^p$ with coefficients in $\cS(\R^p)$.
With the natural identification $\cS(\R^p)\cong \cE^0$, $d_{\R^p}$ the
exterior derivative, and $\tau=\int_{\R^p}$ we obtain a cycle 
$(\cE,d_{\R^p},\int_{\R^p})$ of degree $p$ over $\cS(\R^p)$. 

By \cite[Cor.~5.3]{Ell:CCO} the character of this cycle gives
an isomorphism $HP^{\,p\text{ mod } 2}(\cS(\R^p))\cong \C$. Furthermore
$HP^{\,p+1\text{ mod } 2}(\cS(\R^p))\cong 0$. Therefore, it is appropriate
to call $(\cE,d_{\R^p},\int_{\R^p})$ the \emph{fundamental cycle} of
$\cS(\R^p)$.
\end{example}

Turning back to the cycle $C=(\Omega,d_\cA,\tau)$ we construct
a cycle of degree $n+p$, the \emph{cup-product} of $(\Omega,d_\cA,\tau)$ 
by the fundamental cycle of $\cS(\R^p)$, as follows:
\begin{equation}
  (\Omega^\bullet\cup \cE^\bullet)^k:=\bigoplus_{i+j=k} \Omega^i\otimes \cE^j.
\end{equation}
Note that with  the natural identification $\Omega^i\otimes\cS(\R^p)=
\cS(\R^p,\Omega^i)$, elements of  $\Omega^i\otimes \cE^j$
correspond to differential $j$-forms on $\R^p$ with coefficients
in $\cS(\R^p,\Omega^i)$. We will adopt this point of view if convenient.

The differential on $(\Omega^\bullet\cup \cE^\bullet)^k$ is defined as
\begin{equation}
    \widetilde d(a\otimes \omega)=d_\cA a \otimes \omega+(-1)^i a\otimes
    d_{\R^p}\omega
  \end{equation}
for $a\in \Omega^i$ and $\omega\in \cE^j$. Respectively, if 
$f\in\cS(\R^p,\Omega^i)$, then
\begin{equation}
    \widetilde d(fd\mu_I)= (d_\cA\circ f) d\mu_I+(-1)^i d_{\R^p}(fd\mu_I).
  \end{equation}
Finally, we define a continuous linear functional $\widetilde\tau$
on $(\Omega^\bullet\cup \cE^\bullet)^{n+p}$ by putting
\begin{equation}
   \widetilde\tau(f d\mu_1\wedge\ldots\wedge d\mu_p)=\int_{\R^p} 
   \tau(f(\mu))d\mu
\end{equation}
and $\widetilde\tau_{|\Omega^i\otimes \cE^j}=0$ for $j<p$.
One then checks immediately
\begin{lemma} 
 The cup-product $(\Omega,d_\cA,\tau)\cup (\cE^\bullet,d,\int_{\R^p}):=   
 (\Omega^\bullet\cup\cE^\bullet,\widetilde d,\widetilde\tau)$ 
 is a cycle of degree $n+p$ over $\cS(\R^p,\cA)$.
 The fundamental cycle of $\cS(\R^p)$ is the $p$-fold cup product of the 
 fundamental cycle of $\cS(\R)$ by itself.
\end{lemma}
From \cite[3.3 and 3.7]{Ell:CCO} we now infer
\begin{proposition}
The Elliott--Natsume--Nest-map
\begin{equation}
  \underbrace{\#\circ\ldots\circ\#}_{p \text{ \rm times}}  :  \:
  Z^n_\lambda(\cA)\longrightarrow Z_\lambda^{n+p}(\cS(\R^p,\cA))
\end{equation}
is given by assigning to the character of a cycle $C$ of degree $n$
over $\cA$ the character of the cycle $C\cup(\cE^\bullet,d,\int_{\R^p})$ 
of degree $n+p$ over $\cS(\R^p,\cA)$.
\end{proposition}

Consider now our relative cycle
$(\Omega,\partial \Omega,r,\varrho,\fTR,\lTR)$ 
(cf. \eqref{ML-G2.8}--\eqref{ML-G2.11}). Restricting it to the
algebra of parameter dependent smoothing operators we obtain 
a cycle of degree $p$, $(\Omega,d, \fTR)$, over 
$\CL^{-\infty}(M,E;\R^p)=\cS(\R^p,\CL^\infty(M,E))$. The $L^2$-trace
gives rise to a canonical cycle of degree $0$, $(\CL^{-\infty}(M,E),\tr)$,
over $\CL^{-\infty}(M,E)$ and the cycle $(\Omega,d, \fTR)$ is just
the cup product of $(\CL^{-\infty}(M,E),\tr)$ by the fundamental cycle
over $\cS(\R^p)$. The interesting fact we showed is that the cycle
\begin{equation}
(\Omega,d, \fTR)=(\CL^{-\infty}(M,E),\tr)\cup (\cE^\bullet,d_{\R^p},\int_{R^p})
\end{equation}
\emph{extends} through the exact sequence \eqref{Eq:ExSeqSusAlg} to a relative
cycle $(\Omega_p, \, \partial \Omega_p,\, r,\, \varrho,\, \fTR,\, \lTR)$ 
over the full algebra of
parameter dependent pseudodifferential operators. This was possible since the 
integrated $L^2$-trace $\displaystyle \tr\cup\int_{\R^p}$ has a tracial 
extension, namely the
regularized trace $\fTR$. 
\section{Spectral flow as divisor flow}
As Melrose pointed out from its very inception (cf.~\cite[p. 543]{Mel:EIF}),
the divisor flow has properties which closely parallel those of the spectral flow.
The goal of this last section is to show that this analogy can actually be upgraded
to a precise relationship, which moreover makes sense in every dimension, regardless
of parity.

Since at the $K$-theoretical level the distinction between ``even'' and ``odd'' is 
encoded by the Clifford algebra, we preface the discussion by briefly recalling 
some basic facts about it, which will also allow us to establish the notation.

First, we recall that the Clifford algebra $\C\ell_p$ is the 
universal $C^*$-algebra generated by $p$ unitaries
$e_1,\ldots,e_p$ subject to the relations
\begin{displaymath}
    e_i \cdot e_j + e_j \cdot e_i = -2 \delta_{ij}.
\end{displaymath}
For $p = 2k+1$, $\C\ell_p$ has a unique  irreducible
representation, which in standard form is realized by a $\ast$-homomorphism
 $c: \C\ell_{2k+1} \to \mathfrak{M}_{2^{k}}(\C)$
satisfying
\begin{displaymath}
c ( i^{k+1}\,  e_1 \, \cdots \, e_{2k+1}) \, = \, \id .
\end{displaymath}
When $p =2k$ the standard Clifford representation  $c: \C\ell_{2k} \to \mathfrak{M}_{2^{k}}(\C)$
sends the volume element into a grading operator,
\begin{equation*}
c ( i^{k}\,  e_1 \, \cdots \, e_{2k}) \, = \, \gamma \, , \qquad \text{with} \quad
\gamma^2 = \id \quad \text{and} \quad \gamma^\ast = \gamma \, ,
  \end{equation*}
which gives a decomposition   
\begin{equation}
\label{Eq:DecCliff}
  \C^{2^k} = \Delta^+ \oplus \Delta^-, \quad 
  \text{where } \Delta^\pm = \Ker ( \gamma \mp \id) \, ,
\end{equation}
such that
\begin{equation*}
   c (\mu) = \left(
   \begin{array}{cc}
     0 & -c^+ (\mu)^* \\
     c^+ (\mu) & 0 
   \end{array}
   \right)  \, , \qquad \text{for } \, \mu \in \R^p \, .
\end{equation*}
Here we have been making the identification
\begin{equation}
     \mu \, \equiv \, 
  \sum_{j=1}^p \mu_j \, e_j    \, , \qquad \text{for } \, \mu \in \R^p  \, ,
\end{equation}
which will remain in use for all dimensions $p \in \N$, regardless of parity.

\subsection{Odd case}
Starting with the odd case, we
 consider a smooth path of first order self-adjoint elliptic 
\emph{differential} operators acting between sections 
of $E$, $(D_s)_{0\le s\le 1}$, with 
 $D_0$, $D_1$ invertible, and define its \emph{$p$-fold suspension} as
\begin{equation*}  
  \cD_{p,s}^{\pm}(\mu) :=  D_s \otimes I_{\C^{2^k}}\pm c(\mu)   \, , \qquad \, \mu \in \R^p \, ,
  \quad p = 2k+1 \, .
    \end{equation*}

\begin{theorem} \label{sfodd}
The suspended family 
$\big(\cD_{p,s}^{\pm} \big)_{0\le s\le 1}$
is a smooth path of elliptic elements
in $\CL^1_{p}$, with $\cD_{p,0}$, $\cD_{p,1}$ invertible,
and its divisor flow $\DF$ is related to the spectral flow $\SF$ of the original family by
the identity
\begin{equation} \label{sfdf}
   \DF \bigl( ( \cD_{p,s}^\pm )_{0\le s\le 1} \bigr) =
   \pm \SF\bigl((D_s)_{0\le s\le 1}\bigr) \, .
\end{equation}
\end{theorem}

\begin{proof}  
Since $D_s$ is a \emph{differential operator},
the parametric complete symbol
of $D_s \otimes I_{\C^{2^k}}\pm c(\mu)$ is a polynomial in the cotangent variables and $\mu$,
and therefore $\cD_{p,s}^{\pm} \in\CL^1_{2k+1}$.

If $\sigma_{D_s}(x,\xi)$ denotes the leading symbol of $D_s$ then
$\sigma_{D_s}(x,\xi)\otimes I_{\C^{2^k}} \pm c(\mu)$ is the parametric leading symbol
of $\cD_{p,s}$. Thus $\cD_{p,s}^\pm$ is elliptic since 
$D_s$ is elliptic. Also, we infer from
\begin{equation}
    ( \cD_{p,s}^\pm)^*\cD_{p,s}^\pm \, = \, \Bigl(D_s^2+|\mu|^2\Bigr)\otimes I_{\C^{2^k}}
\end{equation}
that $D_s$ is invertible if and only if $\cD_{p,s}$ is invertible
in $\CL^1_{2k+1}$.
 
The identity Eq.~\eqref{sfdf} now follows by considering the \emph{spectral
$\eta$-invariant} of $D_s$ (cf.~\cite[Sec.~1.13 and 3.8]{Gil:ITHEASIT}). 
Recall that the $\eta$-function of $D_s$,
\begin{equation}
\label{Eq:SpecEtaFunc}
   \eta(z;D_s):=\sum_{\lambda\in\operatorname{spec} D_s} 
   (\operatorname{sgn} \lambda) |\lambda|^{-z},
\end{equation}
extends meromorphically to $\C$ with isolated simple poles. Furthermore,
$0$ is not a pole and one defines the reduced $\eta$-invariant of $D_s$ as
\begin{equation}
     \widetilde{\eta}(D_s):=\frac 12 \bigl(\eta(0;D_s)+\dim\Ker D_s\bigr).
\end{equation}
As a function of $s$, the reduced $\eta$-invariant may have integer jumps.
Hence, $\mod \Z$, the reduced $\eta$-invariant depends smoothly on $s$.
The net number of the integer jumps equals the {\it spectral flow}
\cite[Lemma 3.4]{KirLes:EIM}. Namely,
\begin{equation}
  \SF\bigl((D_s)_{0\le s\le 1}\bigr)=
  \widetilde\eta(D_1)-\widetilde \eta (D_0)-
  \int_0^1\frac{d}{ds}(\widetilde \eta(D_s)\,\textup{mod}\, \Z)\, ds.
\end{equation}

For the moment, consider $s\in J$ in an open subinterval of $[0,1]$, where
$D_s$ is invertible. Then by \cite[Prop. 6.6]{LesPfl:TAP} we have
for $s\in J$
\begin{equation}\label{ML-G2.26}
  \widetilde\eta(D_s)= \pm \frac{k!}{(-2\pi i)^{k+1}(2k+1)!}
  \fTR\Bigl(\bigl((\cD_{p,s}^\pm)^{-1} d\cD_{p,s}^\pm\bigr)^{2k+1}\Bigr).
\end{equation}
The right hand side of Eq.~\eqref{ML-G2.26} has been studied extensively 
in \cite{LesPfl:TAP}. Up to a sign, it was called the (parametric) 
$\eta$-invariant 
of the family $\cD_{p,s}^\pm $. The variation formula for the parametric 
$\eta$-invariant \cite[Prop. 6.3]{LesPfl:TAP}, which also follows from the
transgression formula  Eq.~\eqref{Eq:ChTransgress} and Eqs.~\eqref{Eq:FiEq}, \eqref{Eq:SecEq},
 yields for $s\in J$
\begin{equation}\label{ML-G2.27}\begin{split}
    &\frac{d}{ds} \big(\widetilde \eta(D_s)\,\textup{mod}\,\Z\big)\\
    &= \pm 
    \frac{k!}{(-2\pi i)^{k+1}(2k)!} 
    \lTR\,\Big(\! \sigma(\cD_{p,s}^\pm)^{-1}\sigma(\partial_s 
     \cD_{p,s}^\pm )
    \big(\sigma(\cD_{p,s}^\pm)^{-1}d\sigma(\cD_{p,s}^\pm)\big)^{2k}\Big),
  \end{split}
\end{equation}
where $\sigma$ is the symbol map from the exact sequence \eqref{Eq:ExSeqSusAlg}.
We have used the fact that $\lTR$ is symbolic, cf. Eq.~\eqref{ML-G2.10}.
Hence Eq.~\eqref{ML-G2.27} makes sense if $D_s$ is just elliptic.
A priori, however, Eq.~\eqref{ML-G2.27} does only hold for $D_s$ invertible. 
If $D_s$ is
invertible except for finitely many $s\in [0,1]$ then \eqref{ML-G2.27}
does hold on $[0,1]$ for continuity reasons. The general case is treated
by the usual general position argument:
indeed there exists a sequence $\varepsilon_j>0$, $\varepsilon_j\to 0$ such
that $D_s+\varepsilon_j$ is invertible except for finitely many values of $s$.
Hence \eqref{ML-G2.27} holds for $D_s+\varepsilon_j$ and for continuity 
reasons it does hold for $D_s$ and all $s\in [0,1]$.

Inserting Eqs.~\eqref{ML-G2.26} and \eqref{ML-G2.27} into the formula \eqref{Eq:dfodd}
for the divisor flow (cf. also Eq.~\eqref{df}) gives
\begin{equation}\begin{split}
   \DF&\bigl((\cD_{p,s}^\pm)_{0\le s \le 1}\bigr)\\
   = & \, \frac{k!}{(-2\pi i)^{k+1}(2k+1)!}
   \Big(\fTR \big( \big( (\cD_{p,1}^\pm)^{-1} d\cD_{p,1}^\pm
   \big)^{2k+1}\big)- \fTR \big( \big( (\cD_{p,0}^\pm)^{-1} 
   d\cD_{p,0}^\pm \big)^{2k+1}\big) \Big) \\ &\quad
   -  \frac{k!}{(-2\pi i)^{k+1}(2k)!} 
   \int_0^1 \lTR \left( \big( (\cD_{p,s}^\pm)^{-1} \, 
   \partial_s \cD_{p,s}^\pm \big)\big( (\cD_{p,s}^\pm)^{-1} \, 
   d \cD_{p,s}^\pm\big)^{2k}  \right) \, ds,\\
   =& \pm \left( \widetilde\eta(D_1)-\widetilde\eta(D_0)-\int_0^1\frac{d}{ds} 
   (\widetilde \eta(D_s)\,\textup{mod}\, \Z)\, ds \right)\\
   =& \pm \SF\bigl((D_s)_{0\le s\le 1}\bigr).\qedhere
\end{split}
\end{equation}
\end{proof}

\subsection{Even case}
We begin by considering a single 
operator $D:\cC^\infty (E) \rightarrow \cC^\infty (E)$, 
assumed to be a
first order invertible self-adjoint differential operator acting
on a vector bundle $E\rightarrow M$. 
Its spectral $\eta$-function then satisfies
\begin{equation}
   \eta (s;D) =  \tr \big( D (D^2)^{-(s+1)/2} \big) .
\end{equation}
In \cite[Prop.~6.5]{LesPfl:TAP}, it has been shown for odd $p$ that
\begin{equation}
\label{Eq:EtaRepTr}
  \eta (D) := \eta (0;D) = 
  \frac{\Gamma \big( \frac{p+1}{2}\big)}{\pi^{(p+1)/2}} \,
  \fTR \big(D( D^2 + |\id_{\R^p}|^2)^{-(p+1)/2}\big) ,
\end{equation}
where for reasons of clarity the symbol-valued trace on $\CL (M,E;\R^p)$
has been denoted with a subscript, i.e.~by $\TR$.
In the following we will prove that this formula also holds true for 
even $p$. To this end we first have to recall some analytic tools, 
cf.~\cite[Sec.~6]{LesPfl:TAP}. 

If $f:(0,\infty )\rightarrow \C$ denotes a locally integrable function 
with log-poly\-homogeneous asymptotic expansions for $x\rightarrow 0$ and 
$x\rightarrow \infty$, we put
\begin{equation}
  \regint_0^\infty f(r) dr := \LIM\limits_{\varepsilon \rightarrow 0} 
  \int_\varepsilon^1 f(r)dr + \LIM\limits_{R\rightarrow \infty} \int_1^R f(r) dr,
\end{equation} 
where $\LIM$ stands for the constant term in the corresponding asymptotic 
expansion. Using this regularized integral the following formula has been 
shown in \cite[Prop.~6.5]{LesPfl:TAP} for $z\in \N^*$:
\begin{displaymath}
\begin{split}
  &{z-1-\frac{s+1}{2} \choose z-1 } \eta  (s;D) = \\
   & \! \!=  2  
  \, \frac{\sin \pi\frac{s+1}{2}}{\pi} \regint_0^\infty r^{2z-2-s} 
  \TR_1 \big(D (D^2 +  |\id_{\R}|^2)^{- \, z} \big) (r) dr.
\end{split}
\end{displaymath}
Note that both sides are meromorphic in $s\in \C$. 
Expressing the binomial coefficient and $\frac{\pi}{\sin \pi y}=\Gamma (y) \Gamma (1-y)$
in terms of $\Gamma$-functions, one obtains
\begin{align}
\label{Eq:EtaFunRegInt}
  \eta &(s;D)= \\
       & \! \!= \frac{2\Gamma (z)}{\Gamma \big(\frac{s+1}{2}\big)
  \Gamma (z - \frac{s+1}{2}\big)} 
  \regint_0^\infty r^{2z-2-s} 
  \TR_1 \big(D (D^2 +  |\id_{\R}|^2)^{- \, z} \big) (r) dr.\nonumber
\end{align}
By the argument for the proof of \cite[Prop.~6.5]{LesPfl:TAP}, it is clear that 
this formula actually holds for all real $z>\frac 12$ up to a discrete set
(and that the right side actually extends to a meromorphic function in $z$).
Now observe that under the assumption $z>\frac 12$ the left side is regular at 
$s=0$ and that the factor 
$\frac{\Gamma (z)}{\Gamma \big(\frac{s+1}{2}\big)
  \Gamma (z - \frac{s+1}{2}\big)}  $ is both regular at $s=0$
and non-vanishing. Hence the regularized integral in Eq.~\eqref{Eq:EtaFunRegInt} 
is regular at $s=0$ as well, and we have for $z>\frac 12$, 
\begin{equation}
\label{Eq:EtaRegInt}
  \eta (D) = \frac{2\Gamma (z)}{\sqrt{\pi} \, 
  \Gamma (z - \frac{1}{2}\big)} 
  \regint_0^\infty r^{2z-2} 
  \TR_1 \big(D (D^2 +  |\id_{\R}|^2)^{- \, z} \big) (r) dr.
\end{equation}
 Using this, and the rotation invariance of the function (defined on $\R^p$)
 $\fTR\big( D (D^2 +|\id_{\R^p}|^2)^{- \frac{p+1}{2}} \big)$,
 one obtains:
\begin{align}
  \fTR & \, \big( D (D^2 +|\id_{\R^p}|^2)^{- \frac{p+1}{2}} \big)\nonumber\\
  & =
  \regint_{\R^p} \TR \big( D (D^2 + |\id_{\R^{p}}|^2)^{- \frac{p+1}{2}} \big) 
  (\mu) \, d\mu \label{Eq:etatrbar} \\
  & = \LIM\limits_{R\rightarrow \infty} 
  \frac{p\,\pi^{p/2}}{\Gamma\big( \frac p2 +1\big)}
  \int_0^R r^{p-1} \TR_1 \big( D (D^2 + |\id_{\R}|^2 )^{- \frac{p+1}{2}} \big) 
  (r) \, dr \nonumber\\
  & = \frac{p\,\pi^{p/2}}{\Gamma\big( \frac p2 +1\big)}
  \regint_0^\infty r^{p-1} 
  \TR_1 \big( D (D^2 + |\id_{\R}|^2 )^{- \frac{p+1}{2}}\big) (r) \, dr \nonumber\\
  & = \frac{\pi^{(p+1)/2}}{\Gamma \big( \frac{p+1}{2} \big)} \, \eta (D) .\nonumber
\end{align}
This shows that Eq.~\eqref{Eq:EtaRepTr} holds in all dimensions $p$.

We now use the standard Clifford representation $c$ to
define the $p$-fold suspension of 
the operator $D$ in the even case $p=2k$ as the 
parametric differential operator 
\begin{equation}\label{Eq:Dsuspensioneven}
  \cD_{2k}(\mu) : =  \gamma \, \big( D  \otimes I_{\C^{2^k}} + c(\mu) \big) 
  \,  = \, 
  \left( \begin{array}{cc}
     D & -c^+ (\mu)^* \\
     - c^+ (\mu) &  -D 
   \end{array}
   \right) .
\end{equation}
By construction, $\cD_{2k}$ is an element of $\CL^1 (M,E \otimes \C^{2^k};\R^{2k})$.
From the invertibility of $D$, it follows that the operator $\cD_{2k}$ is 
invertible. Moreover, $\cD_{2k}(\mu)^2 = \bigl(D^2 + |\mu|^2\bigr)\otimes  I_{\C^{2^k}}$ is diagonal 
with respect to the decomposition \eqref{Eq:DecCliff}. Hence 
\begin{displaymath}
  \cQ := \big( \cD_{2k} ^2\big)^{1/2} = \big(D^2 + | \id_{\R^{2k}}|^2\big)^{1/2}
  \otimes I_{\C^{2^k}} \in \CL^1 (M,E\otimes \C^{2^k};\R^{2k}) 
\end{displaymath}
is invertible and
\begin{equation}\label{Eq:Dsusid}
  \cP := \frac 12 \big( I - \cQ^{-1}\cD_{2k} \big) \in 
  \CL^0 (M,E\otimes \C^{2^k};\R^{2k}) 
\end{equation}
is an idempotent. 

Let us now determine $\fTR \Big( \big( \cP-\frac 12\big) (d\cP)^{2k}\Big)$. 
To this end observe first that $\cD_{2k}$ commutes with $\cQ$ and that
\begin{equation}
\label{Eq:DerOpcD}
  d\cD_{2k} = \sum_{j=1}^{2k} \gamma \, c(e_j) \, d \mu_j .
\end{equation}
One also checks immediately 
\begin{displaymath}
  d\cQ \wedge d\cQ =0, \quad \cQ^{-1} d\cD_{2k} = (d\cD_{2k}) \cQ^{-1} , \quad 
  d \cQ \wedge d\cD_{2k} + d \cD_{2k} \wedge d \cQ =0 .
\end{displaymath}
These relations entail the following two chains of equalities:
\begin{displaymath}
\begin{split}
  d\cP \wedge d\cP  = & \,\frac 14 d(\cQ^{-1}\cD_{2k} )\wedge d (\cQ^{-1}\cD_{2k}) = \\
  = \, &\frac 14 \big( \cQ^{-2} (d\cQ) \cD_{2k} \wedge \cD_{2k} \cQ^{-2} (d \cQ) - 
  \cQ^{-2} (d\cQ) \cD_{2k} \wedge \cQ^{-1} d \cD_{2k} - \\
  & - \cQ^{-1} d \cD_{2k} \wedge \cD_{2k} \cQ^{-2} (d\cQ)  + \cQ^{-1} d\cD_{2k} \wedge 
  \cQ^{-1} d\cD_{2k} \big) \\
  = \, &\frac 14  \cD_{2k}^{-2} d \cD_{2k}  \wedge d\cD_{2k} ,
\end{split}
\end{displaymath}
and
\begin{displaymath}
\begin{split}
  (d\cP)^{2k} \, & =  4^{-k} \, \cD_{2k}^{-2k} \, (-1)^k \,
  \sum_{\sigma \in S_{2k}} \, c(e_{\sigma (1)})\,
  \cdots \, c(e_{\sigma(2k)}) \, d\mu_{\sigma (1)} \wedge \ldots \wedge
  d\mu_{\sigma(2k)} \\
  & =  (2k)! \, 4^{-k} \, \cD_{2k}^{-2k} \, i^k \, \gamma \,
  d\mu_1 \wedge \ldots \wedge d\mu_{2k} .
\end{split}
\end{displaymath}
Hence we obtain
\begin{equation}
\begin{split}
  \big(\cP -\frac 12\big) (d\cP)^{2k} \, &=  \frac{ - (2k)!\, i^k }{2^{2k+1}}
  \cQ^{-1} \cD_{2k}^{-2k}\gamma\big( D +c(\mu)\big)\gamma \, 
  d\mu_1\wedge\ldots\wedge d\mu_{2k}\\
  & = \frac{-(2k)! \, i^k}{2^{2k+1}}\cQ^{-2k-1}\big( D-c(\mu)\big)
  \, d\mu_1 \wedge \ldots \wedge d\mu_{2k} ,
\end{split}
\end{equation}
and finally
\begin{equation}\label{Eq:etaequalspareta}
\begin{split}
  \fTR&\Big( \big( \cP -\frac 12 \big) (d\cP)^{2k} \Big) \\
       &= \frac{-(2k)! \, 
  i^k}{2^{k+1}}
  \regint_{\R^{2k}} \TR \big( D (D^2 +|\mu|^2)^{-(2k +1)/2} \big) \, d\mu \\
  & = -\frac 12 (2\pi i)^k k! \, \eta (D).
\end{split}
\end{equation}
Here we have used \eqref{Eq:etatrbar} and the fact that the standard
Clifford representation has rank $2^k$.

Summing up we have proved the following even analogue of
\cite[Prop. 6.6]{LesPfl:TAP}.

\begin{proposition}\label{etaequalspareta}
Let $D$ be an invertible first order self-adjoint elliptic differential
operator. Let $\cD_{2k}(\mu):= \gamma \, \big( D  \otimes I_{\C^{2^k}} + c(\mu) \big), 
\mu\in \R^{2k},$ be the $2k$-fold
suspension defined in Eq.~\eqref{Eq:Dsuspensioneven} and let $\cP$ be
the idempotent defined in Eq.~\eqref{Eq:Dsusid}. Then the $\eta$-invariant
of $D$ satisfies
\begin{displaymath}
  \eta(D)\, =\, -\frac{2}{ (2\pi i)^k k! }
 \fTR\Big( \big( \cP -\frac 12 \big) (d\cP)^{2k} \Big).
\end{displaymath}
\end{proposition}

\begin{remark} This identity justifies promoting  the above expression to a definition:
  the {\it higher (even) $\eta$-invariant}  $ \eta_{2k}$ is defined on 
 projections $\cP\in \CL^0 (M,E;\R^{2k})$ by
 \begin{equation}
   \eta_{2k} (\cP) := - \frac{2}{(2\pi i)^k k!} \fTR \,  
   \Big( \big( \cP -\frac 12 \big) (d\cP)^{2k}  \Big) .
 \end{equation}
\end{remark}
%


Finally, we record the even analogue of Theorem \plref{sfodd}. 

\begin{theorem} \label{sfeven}
 Let $D_s :\cC^\infty (E) \rightarrow \cC^\infty (E)$ be a smooth family of elliptic
 first order self-adjoint differential operators on the vector bundle $E$
 such that $D_0$ and $D_1$ are invertible.
Let $\cD_s:=\gamma \, \big( D_s  \otimes I_{\C^{2^k}} + c(\mu) \big), \mu\in\R^{2k},$
be the $2k$-fold suspension defined in Eq.~\eqref{Eq:Dsuspensioneven}.
Furthermore, let $\cP_s\in \CL^1 (M,E\otimes \C^{2^k};\R^{2k}) $,
$s \in [0,1]$,  be a smooth
family of almost idempotents with endpoints
\[
       \cP_j=\frac 12\bigl( I- (\cD_j^2)^{-1/2}\cD_j\bigr),\quad j=0,1
\]
and whose symbols satisfy
\[
    \sigma(\cP_s)=\frac 12 \bigl(I- \sigma(\cD_s^2)^{-1/2}\sigma (\cD_s)\bigr).
\]
Then the even divisor flow of the family of almost idempotents 
 $(\cP_s)_{s \in [0,1]}$ coincides with the spectral flow of $(D_s)_{s \in [0,1]}$:
 \begin{equation}
   \DF \big( (\cP_s)_{0\leq s\leq 1}\big) = 
   \SF \big( (D_s)_{0\leq s\leq 1}\big).
 \end{equation}
\end{theorem}

\begin{proof} We first prove the existence of a family of almost idempotents
$(\cP_s)_{0\leq s\leq 1}$ with the stated properties. Start with the smooth family
\begin{equation}
    p_s=\frac 12 \bigl(I- \sigma(\cD_s^2)^{-1/2}\sigma (\cD_s)\bigr)
\end{equation}
of projections in $\csym^0_{2k}$. In view of Eq.~\eqref{Eq:symbolsection} and
Eq.~\eqref{Eq:symbiso} there is a smooth lift $\widetilde \cP_s \in \CL^1 (M,E\otimes \C^{2^k};\R^{2k}) $
with $\sigma(\widetilde \cP_s)=p_s$. To adjust the endpoints we put
\begin{equation}\label{Eq:adjendpts}
   \cP_s:=\widetilde \cP_s+s (\cP_1-\widetilde\cP_1)+(1-s)(\cP_0-\widetilde \cP_0),
\end{equation}
which has all the desired properties. 

Alternatively and even more concretely one can obtain $\cP_s$ by modifying the
construction of $\cQ_s$ as follows: choose an even smooth function $\phi\in \cC_\text{\tiny c}^\infty(\R)$
with compact support such that $\phi(0)=1$ and put
\begin{equation}
    \cQ_s:= (\cD_s^2+\phi(\mu)\phi(D_s^2))^{1/2}.
\end{equation}
$\cQ_s$ is invertible and Eq.~\eqref{Eq:Dsusid} now yields a smooth family of almost idempotents
$\widetilde \cP_s$ with $\sigma(\cP_s)=p_s$. The endpoints are adjusted as before in
Eq.~\eqref{Eq:adjendpts}.


Once such a family of almost idempotents is chosen, 
the proof of the statement is completely analogous to that
of Theorem \plref{sfodd}. The only difference is that
for the variation of the even parametric $\eta$--invariant we cannot refer
to \cite{LesPfl:TAP} but have to use the transgression formula
Eq.~\eqref{evenh} and Eqs.~\eqref{Eq:FiEqeven}, \eqref{Eq:SecEqeven}.
More precisely, assume for the moment that $D_s$ is invertible and
let $(\varphi_{2k},\psi_{2k+1})$ be the character
of the cycle $C^{2k}_\text{\rm\tiny reg}$ (cf. Eq.~\eqref{regcycle}). Then from Proposition \ref{etaequalspareta}
and Eq.~\eqref{Eq:FiEqeven} we infer
\begin{equation}\label{ML-G2.26even}
\begin{split}
  \widetilde\eta(D_s)\, &=\, -\frac{1}{ (2\pi i)^k k! } \fTR\Big( \big( \cP_s -\frac 12 \big) (d\cP_s)^{2k} \Big),\\ 
         &=\frac{(-1)^{k+1}}{ (2\pi i)^k} \langle \varphi_{2k},\ch_\bullet(\cP_s)\rangle
\end{split}
\end{equation}
and hence Eq.~\eqref{evenh}, Proposition \ref{ML-prop3.1} and Eq.~\eqref{Eq:SecEqeven}
yield
\begin{equation}\label{ML-G2.27even}\begin{split}
    &\frac{d}{ds} \big(\widetilde \eta(D_s)\,\textup{mod}\,\Z\big)\\
    &= \frac{(-1)^{k+1}}{ (2\pi i)^k}\langle \varphi_{2k},(b+B)\slch_\bullet(\cP_s,(2\cP_s-1)\partial_s\cP_s))\rangle\\
    &= \frac{(-1)^{k+1}}{ (2\pi i)^k}\langle \sigma^*\psi_{2k-1},\slch_\bullet(\cP_s,(2\cP_s-1)\partial_s\cP_s))\rangle\\
    &=\frac{-1}{ (2\pi i)^k k!}\lTR\Bigl(\sigma(2\cP_s-1)\sigma(\partial_s\cP_s)(d\sigma(\cP_s))^{2k-1}\Bigr).
  \end{split}
\end{equation}
Thus we have established the analogue of Eq.~\eqref{ML-G2.27} in the even case. 
As in the proof of Theorem \ref{sfodd} one now shows that the equality between the first
and the last term in Eq.~\eqref{ML-G2.27even} holds
for elliptic families $D_s$. 
Then inserting Eqs.~\eqref{ML-G2.26even} and \eqref{ML-G2.27even} into the formula \eqref{Eq:dfeven}
for the divisor flow gives the claim.
\end{proof}

%

%

\end{document}